\documentstyle{article}
\begin{document}
\footnotetext{Supported in part by the Bulgarian
foundation ``Scientific Research'' under the
project MM-404/94}

\centerline{\sc  MINIMAL SECTIONS OF CONIC BUNDLES}

\vspace{0.8cm}

\centerline{\tt ATANAS ILIEV}

\vspace{1.3cm}

\begin{abstract}
%{\footnotesize
Let $X$ be a smooth conic bundle over the projective plane,
and let $(J(X),{\Theta})$ be the principally polarized
intermediate jacobian of $X$.
In this paper we relate the Wirtinger description of
$(J(X),{\Theta})$ as a Prym variety, and
the description of maximal subbundles of rank 2 vector bundles
over curves, to prove the existence of two canonical families
${\cal C}_{+}$ and ${\cal C}_{-}$ of curves on $X$ such that
the Abel-Jacobi map
sends one of these families onto a copy of the theta divisor
${\Theta} \subset J(X)$, and the other -- onto
$J(X)$. The elements of ${\cal C}_{+}$
and ${\cal C}_{-}$ correspond to the non-isolated and to the
isolated
minimal sections of a naturally defined system of conic bundle
surfaces on $X$. In either of the two possible
cases we describe the
general fiber of the Abel-Jacobi map for ${\cal C}_{+}$ and for
${\cal C}_{-}$.
As an application, we study the families ${\cal C}_{+}$ and
${\cal C}_{-}$
and their Abel-Jacobi images in the setting of
the natural birational conic bundle structures on some
classes of Fano threefolds.
In particular, we find special parameterizations of the
theta divisor of the cubic threefold, of the bidegree (2,2)
threefold, and of the nodal quartic double solid (q.d.s.).
The last involves the Tikhomirov's parameterization of
the theta divisor for the general q.d.s. $X$:
$\Theta = {\Phi}_{\cal R}({\cal R})$, where  ${\cal R}$ is the
family of Reye sextics on $X$.
In addition, the theory of minimal sections gives that
the general fiber of ${\Phi}_{\cal R}$ has exactly two
connected components.
%}
\end{abstract}

\bigskip
\bigskip

\centerline{\bf Preliminaries}

\smallskip
%
%{\it General remarks}
%
%\smallskip

In the classification of algebraic threefolds,
conic bundles take a special place.
One of the important applications of the Mori theory of
minimal models, due especially to investigations of
Miyaoka (see e.g. [Miy]), is the obtained classification
of regular threefolds of negative Kodaira dimension.
In general, any such
threefold $V$ is birational to a threefold $X$
with at most terminal singularities and of (at least) one
of the following types: 1. {\bf Q}--Fano threefolds;
2. Conic bundles $X \rightarrow S$ over normal surfaces $S$;
3. Del-Pezzo bundles $X \rightarrow C$ over smooth curves $C$
(see [Isk] for more details).

Especially, if $X \rightarrow S$ is a conic bundle, and if the
base surface $S$ is not rational, then $X$ is not rational.
Therefore, from the point of view of the problem of
rationality, it makes sense to study conic bundles with a
rational base $S$. Since any rational surface $S$ is birational
to the projective plane ${\bf P}^{2}$, letting
$S = {\bf P}^{2}$ is not a substantial restriction.

{}From some other point of view, any smooth algebraic threefold
$Y$ admitting a rational map with connected and rational
general fibers $f: Y \rightarrow T$
to a rational surface $T$, is birational to a {\it standard}
conic bundle $p:X \rightarrow S$ over a rational surface $S$
(see [Z], [S]).

Let $X$ be a smooth projective threefold, and let $S$ be
a smooth surface.  By definiton, $X$ is a
{\it standard conic bundle} over $S$, if there exists
a surjective morphism $p:X \rightarrow S$ such that:
(i).The general fiber ${p}^{-1}(s)$ is
a smooth connected rational curve (a conic);
(i).For any irreducible
curve $C \subset S$, the preimage ${p}^{-1}(C) \subset X$
is irreducible.

Let $X$ be standard, let
$\Delta = \{ s \in S : {p}^{-1}(s) \mbox{ is singular} \}$
be the {\it discriminant} set of $p$, and let
$\Delta \neq \oslash$. Then (see e.g. [B1]) $\Delta$ is a
curve with at most double points, and if $s \in {\Delta}$
then $p^{-1}(s)$ is either: (i). a union of two smooth rational
curves (lines) intersecting each other simply in a uniuque
point, and $s$ is a smooth
point of $\Delta$; or  (ii). $p^{-1}(s)$ is a rational curve
counted twice (a double line), and $s$ is a double point
of $\Delta$.

As foolows from the preceding, it makes sense to study
the standard conic bundles over rational surfaces $S$
(esp. $S = {\bf P}^{2}$). Moreover, standard argument
regarding general positions gives that the case
$Sing({\Delta}) \neq \oslash$ can be regarded as
a ``degeneration of the general case''
$\Delta - \mbox{ non-singular}$.

Taking in mind all this, we shall study in detail
the standard conic bundles $p:X \rightarrow {\bf P}^{2}$
with smooth discriminant curves $\Delta$.

Let $p: X \rightarrow {\bf P}^{2}$ be such a conic bundle.
Being a rational fibration over a rational surface,
$X$ is a threefold with a non-effective canonical class,
i.e. $h^{3,0}(X) = h^{0}(X, {\Omega}^{3}_{X}) = 0$, i.e.
the complex torus (the Griffiths intermediate jacobian )
$J(X)$ of $X$ does not contain a $(3,0)$-part. In particular

$J(X) = H^{2,1}(X)^{*}/(H_{3}(X,{\bf Z}):mod.{\sl torsion})$
is a
principally polarized abelian variety (p.p.a.v.)
with a principal polarization (p.p.)
defined by the
intersection of real 3-chains on $X$ (see [CG]).
The divisor $\Theta$ of this polarization
is called the theta divisor of $J(X)$.

Since $p:X \rightarrow {\bf P}^{2}$ is standard and
$\Delta$ is smooth, the splitting
$p^{-1}(s) = {\bf P}^{1} \bigvee {\bf P}^{1} , s \in {\Delta}$
defines a unbranched double covering
${\pi}: \tilde{\Delta} \rightarrow \Delta$ of the
smooth discriminant curve $\Delta$.
Therefore
the pair $(\tilde{\Delta}, {\Delta})$
defines in a natural way the p.p.a.v.
$P(\tilde{\Delta},{\Delta})$ -- the Prym variety of
${\pi}:{\tilde{\Delta}} \rightarrow {\Delta}$
(see section 4);
and the well-known result of Beauville [B1]
tells that $(J(X), {\Theta})$ and
$P(\tilde{\Delta},{\Delta})$ are isomorphic as
p.p.a.v.

The Beauville's approach to conic bundles is
based on the Wirtinger description of
the p.p. Prym variety $P(\tilde{\Delta},{\Delta})$
by sheaves on the double discriminant curve
$\tilde{\Delta}$ -- this way, this approach is
``nearer'' to the associated curve data,
and ``far'' from the geometry of the threefold $X$.
In this paper
we study conic bundles from an alternative
point of view, based on the properties of some special
families of curves on the threefold $X$ -- the
sections of the conic bundle map $p$.
The advantages of this second approach are that it
makes possible to describe the pair
$(J(X),{\Theta})$ in terms of families of curves
on $X$ and their Abel-Jacobi images.

More concretely, let $X$ be a smooth threefold
with $h^{3,0} = 0$, let $(J(X), {\Theta})$ be the
p.p. intermediate jacobian of $X$,
and let $A_{1}(X)$ be the group of rational equivalence
classes of algebraic 1-cycles $C$ on $X$ which are
homologous to $0$.  Then the integrating over the
real 3-chains $\gamma$ s.t.
${\delta}({\gamma})$ = (the boundary of $\gamma$) = $C$,
$C \in A_{1}(X)$ defines the natiural map:
${\Phi}: A_{1}(X) \longrightarrow J(X)$ --
the Abel-Jacobi map for $X$ (see e.g. [CG]).
In addition, if ${\cal C}$ is a smooth family of homologous
cycles $C$ on $X$, and $C_{0}$ is a fixed element
of ${\cal C}$, then the composition of $\Phi$ and the
cycle-class map
$F \rightarrow A_{1}(X)$,
$C \mapsto [C - C_{0}]$,
defines a map ${\Phi}_{\cal C}:{\cal C} \longrightarrow J(X)$.

Let $Alb({\cal C})$ be the Albanese variety of $F$.
By the universal property of the Albanese map
$a: {\cal C} \rightarrow Alb({\cal C})$, ${\Phi}_{{\cal C}}$
can be factorized through $a$, and defines the map
${\Phi}^{'}_{\cal C} : Alb({\cal C}) \longrightarrow J(X)$.

Both ${\Phi}_{\cal C}$ and ${\Phi}^{'}_{\cal C}$ are called
the {\it Abel-Jacobi maps} for the family of
1-cycles ${\cal C}$.

For a large class of such threefolds $X$ (esp. -- for
conic bundles), the transpose  ${}^{t}{\Phi}$ of the
Abel-Jacobi map for $X$ defines an isomorphism between
the Chow group $A_{1}(X)$ and $J(X)$
(see [BM]), and one may expect that
for some ``rich'' families of curves ${\cal C}$ on $X$
the Abel-Jacobi map ${\Phi}_{\cal C}$ will be surjective.
Moreover, one can set the folowing problem:

{\bf (*).}
The problem of parameterization of $\Theta$.

{\it
Find an appropriate family ${\cal C}_{\theta}$ of
algebraic 1-cycles on $X$ such that the Abel-Jacobi
map ${\Phi}_{{\cal C}_{\theta}}$ sends ${\cal C}_{\theta}$
surjectively onto a copy of the theta divisor ${\Theta}$.
}

Assume the existence of such a family ${\cal C}_{\theta}$.
One can
formulate the folowing additional questions:

{\bf (**).}
{\it
Describe, in terms of ${\cal C}_{\theta}$ and $X$, the
structure of the general fiber of ${\Phi}_{{\cal C}_{\theta}}$;
}

{\bf (***).}
{\it
Describe, in terms of ${\cal C}_{\theta}$
and $X$,
(the) components
of the set $Sing({\Theta})$ of singular points of
the theta divisor $\Theta$, and the tangent cones to
$\Theta$ in these singularities.
}

Especially, the answer of (***) is closely related to
the Torelli problem for $X$ (see [CG], [B1], [B2],
[Vo], [De], [I1] etc.).

In this paper we give a positive answer of the problems (*)
and (**) for the general conic bundle
$p:X \rightarrow {\bf P}^{2}$.
More concretely, we prove the existence of two naturally
defined families of connected curves
${\cal C}_{+}$ and ${\cal C}_{-}$ on
the conic bundle $X$, such that
the Abel-Jacobi map sends one of these two families
onto a copy of the theta divisor ${\Theta}$, and
the second -- onto the intermediate jacobian
$J(X)$.
Let $deg ({\Delta}) = d > 3$. It turns out that
these two families of curves are the two components
of the canonical family ${\cal C}_{min}$
of minimal sections of the
conic bundle map $p: X \rightarrow {\bf P}^{2}$ --
the elements of ${\cal C}_{min}$ being the  the minimal
sections of the induced system of ruled surfaces
$\{ p^{-1}(C_{0})$: $C_{0}$ -- a plane curve of degree $d-3 \}.$

The answer of the question:
{\sl
(Which one of
${\cal C}_{+}$ and ${\cal C}_{-}$ parameterizes ${\Theta}$)
}
closely depends on the geometry of the given conic
bundle.  In particular:  the ``theta''-family in Example (6.5)
(the bidegree (2,2) threefold)
is ${\cal C}_{+}$, while the ``theta''-family in Example (6.6)
(the nodal quartic double solid)
is ${\cal C}_{-}$.

Our particular approach to the minimal sections of the conic
bundle $p:X \rightarrow {\bf P}^{2}$ uses the relation
between: 1.the families of sections on ruled surfaces, and  2.the
families of subbundles of rank 2 vector bundles describing
these surfaces; and also -- their deformation theory
(see [LN], [Se]). The general observation is that it is
possible to relate the effective divisors from the
linear systems of the sheaves on $\tilde{\Delta}$
reprezenting $(J(X), {\Theta})$ as a Prym variety, and
the families of minimal sections on some special
non-minimal ruled surfaces $S \subset X$.

As regards the question (***), we give its answer in
various examples: the cubic threefold, the Fano threefold
$X_{16}$, the bidegree (2,2) divisor, the nodal
quartic double solid (see (6.3),(6.4),(6,5),(6.6)).
Although the structure of $Sing({\Theta})$ for these
threefolds is well-known (see e.g. [GH], [Tju],
[B2], [Ve], [I1], [Vo], [De]), the approach
presented here suggests to study these
special subvarieties of $\Theta$ by studying special
subfamilies of degenerate 1-cycles, in the families
of minimal sections ${\cal C}_{+}$ and ${\cal C}_{-}$.

\smallskip

In addition, the presented approach can be easily
generalized to arbitrary standard conic bundles
$X \rightarrow S$ over rational surfaces $S$.
Morerover, the general Prym variety can be
represented (non-uniquely) by the intermediate
jacobian of such a conic bundle
(see [Isk, Lemma 1(iv)]).
Therefore  the minimal sections can be used for
future studying of associated Prym varieties
(see e.g. Example (6.5) -- originated from the
Verra's counterexample to the Torelli theorem for
Prym varieties [Ve] --  where the Prym variety,
of a double covering of a general plane sextic
$\Delta$, is represented by the jacobian of a bidegree
(2,2) threefold $X_{2,2}$; this way -- involving
additional information about $Sing({\Theta})$;
see also [I1]).

\bigskip

%{\bf 0. Introduction.}
%
%\smallskip
%
%{\bf (0.1)}
%{\it General notation.}
%
%\smallskip
%
%$p:X \rightarrow {\bf P}^{2}$ --
%a general conic bundle threefold $X$ over ${\bf P}^{2}$;
%
%$(J(X),{\Theta})$ --
%the p.p. intermediate jacobian of the threefold $X$;
%
%${\pi}: \tilde{\Delta} \rightarrow {\Delta}$ --
%the associated to $p:X \rightarrow {\bf P}^{2}$
%double covering of discriminant curves;
%
%$P(\tilde{\Delta},{\Delta})$ --
%the Prym variety of $(\tilde{\Delta},{\Delta})$;
%
%${\cal C}_{-}$ and ${\cal C}_{+}$ --
%the canonical families
%of isolated and non-isolated minimal sections of $p$
%(see (4.1.2));
%
%${\Phi}_{-}$ and ${\Phi}_{+}$ --
%the Abel-Jacobi maps for the families
%${\cal C}_{-}$ and ${\cal C}_{+}$.

\bigskip

%{\bf (0.2)}
%+++++++++++++++++++++
%
%
\centerline{\bf 0. Structure of the paper.}
\smallskip

Let $S$ be a smooth rational surface.
It is well-known (see e.g. [B1]) that the principally
polarized intermediate jacobian $(J(X),{\Theta})$
of a smooth conic bundle
$p:X \rightarrow S$
is isomorphic to the Prym variety $P$ of the induced
discriminant pair
$(\tilde{\Delta},{\Delta})$.

We shall study in detail the conic bundles
$p:X \rightarrow {\bf P}^{2}$.

We let $S = {\bf P}^{2}$, and
$p:X \rightarrow {\bf P}^{2}$;
$d = deg({\Delta})$,  ${\Delta} \subset {\bf P}^{2}$.

\smallskip

In this paper we prove the existence of two naturally
defined families ${\cal C}_{-}$ and ${\cal C}_{+}$
of homologically equivalent algebraic 1-cycles
(curves) on $X$ (see e.g. (4.1.2))
with the following properties ($d = deg({\Delta}) \ge 4$):

Let
${\Phi}_{-}:{\cal C}_{\epsilon} \rightarrow J(X)$
be, as above,
the Abel-Jacobi map for ${\cal C}_{\epsilon}$
${\epsilon} \in \{ -,+ \}$.
Then one of the folowing two alternatives is true
(see (4.3), (4.4)):

{\bf (A,+)}
${\Phi}_{-}$ is surjective, and ${\Phi}_{+}$ maps
${\cal C}_{+}$ onto a copy of the theta divisor
${\Theta} \subset J(X)$, or

{\bf (A,--)}
${\Phi}_{+}$ is surjective, and ${\Phi}_{-}$ maps
${\cal C}_{+}$ onto a copy of the theta divisor
${\Theta} \subset J(X)$.

{\bf (B)}
The general element
$C \in {\cal C}_{-}$
(resp. -- of ${\cal C}_{+}$)
is maped isomorphically onto a smooth curve
$p(C)$ such that
$p(C) \cap {\Delta}$
is an effective divisor of
the canonical system
$\mid {\omega}_{\Delta} \mid$
of the curve $\Delta$.
Moreover, if
$C \in {\cal C}_{+}$
(resp. -- if $C \in {\cal C}_{+}$)
then $C$ is an isolated
(resp. -- a non-isolated)
minimal section
of the non-minimal ruled surface
$p:{p}^{-1}(p(C)) \rightarrow p(C)$
(see the definitions in (2.2)).

\smallskip

In sections 1 -- 5 we prove the existence,
and describe the general properties, of the
two canonical families of minimal sections
${\cal C}_{-}$ and
${\cal C}_{+}$.

More concretely, let
$p:X \rightarrow {\bf P}^{2}$ be a general
conic bundle, let $\Delta$ be the discriminantal
curve of $p$, and let $deg({\Delta}) = d$.
Let $C$ be e.g. a reduced plane curve of degree $k < d$,
and let $S_{C} = p^{-1}(C)$.
Then $p$ induces a conic bundle structure
$p:S_{C} \rightarrow C$.
This way $p$ defines, for any $k < d$,
the family ${\cal S}[k]$ of conic bundle surfaces
over the base space $\mid {\cal O} (k) \mid$
of plane curves of degree $k$.

Theorem
(2.4) tells that if $X$ is sufficiently general
then the general element $S_{C}$ of
any of the families ${\cal S}[k], k < d$, can
be regarded as a general element of a versal deformation
of a conic bundle surface over a plane curve of degree $k$.
The minimal sections of such a surface $S_{C}$ are of two
types -- isolated and non-isolated -- and correspond to the
minimal sections of the general element $S_{\sigma}$
in any of the two types of a versal deformation of
ruled surfaces
over a plane curve of
degree $k$
(resp. -- of genus $g = (k-1)(k-2)/2$).
These two types are separated by
the parity of the invariant $e(S_{\sigma})$.
The general element $S_{\sigma}$
of the even type: $e \equiv 0 \mbox{ (mod. } g)$
has a $1$-dimensional family of minimal sections,
while the general element of the odd type
has a finite number of minimal sections.
By collecting  all these sections  we find, for any
$k < d$,  two naturally defined families --
${\cal C}_{-}(k)$ and ${\cal C}_{+}(k)$ -- of isolated and
non-isolated minimal sections of the elements of the
family of surfaces ${\cal S}[k]$.

In section 3 we prove Theorem (2.4). The proof is based on:
(a). the versality of the family
${\cal S}[1] \rightarrow$
$\{ \mbox{the lines in } {\bf P}^{2} \}$
(see section 1);
(b). the versality of the families
${\cal S}_{k} \subset {\cal S}[k]$
of conic bundle surfaces over
degenerate curves
$C_{0}$  =
$\mbox{ a union of } k \mbox{ lines in } {\bf P}^{2}$
(see (3.3)).
This involves the versality of ${\cal S}[k]$.

Let ${\cal C}_{-}$ and ${\cal C}_{+}$ be the two families
of minimal sections for $k = d-3$ --
the canonical families of isolated and non-isolated
minimal sections of $p:X \rightarrow {\bf P}^{2}$
(see (4.1.2), the definitions in (2.2), (5.2) -- (5.7)).
Theorem (4.4) tells
that the Abel-Jacobi map sends one of these two
families surjectively onto the intermediate jacobian
$J(X)$, while the second family is mapped onto a
translate of the theta divisor ${\Theta}$.

Theorem (5.8) describes the general fibers of the
Abel-Jacobi maps
${\Phi}_{+}$ and ${\Phi}_{-}$ for ${\cal C}_{+}$ and
${\cal C}_{-}$,
depending on which one of these two families parameterizes
$\Theta$.

\smallskip

In section 6 we study some of the most tipical examples
of conic bundles with $deg({\Delta})$ = 4,5 and 6.

In particular, in example (6.4) we find a new parameterization
of the theta divisor ${\Theta}$ for the cubic threefold
$X_{3}$ -- via the
Abel-Jacobi image of the 6-dimensional family ${\cal C}^{0}_{3}$
of rational cubic curves on $X_{3}$.  Comparing with the known
parameterization of ${\Theta}$ by the Abel-Jacobi image
of the 4-dimensional family
$\{ \mbox{differences of two lines on } X_{3} \}$, the new
parameterization has the advantage that the fiber of the
Abel-Jacobi map for ${\cal C}^{0}_{3}$ is connected, and the
degree of the Gauss map for ${\Theta}$ can be easyly checked to
be 72 = the number of ${\bf P}^{2}$-systems of rational cubics on
a general cubic surface (see (6.3.4)).

In examples (6.4), (6.5) and (6.6), we study the
conic bundle structures on the Fano 3-fold
$X_{16} \subset {\bf P}^{10}$, on the
bidegree (2,2) divisor $X_{2,2}$ in
${\bf P}^{2} \times {\bf P}^{2}$, and on the nodal
quartic double solid.  As an application
we give another proof of the result of
Tikhomirov: {\sl The family ${\cal R}$ of Reye sextics
(sextics of genus 3) on ${X'}_{2}$
parameterizes the theta divisor of the general quartic
double solid ${X'}_{2}.$}
Moreover, we find a natural family of curves which
parameterizes (via the Abel-Jacobi map)
the intermediate jacobian of ${X}'_{2}$
(see (6.6.6) and (6.6.7)(ii)).

It follows from [T] that any connected component of the
general fiber of the Abel-Jacobi map
${\Phi}_{\cal R} : {\cal R} \rightarrow {\Theta}$
is ismorphic to ${\bf P}^{3}$.
Now, the
``minimal section'' approach tells that the number of
these fibers is {\it two}, i.e. the Stein quotient of
${\Phi}_{\cal R}$ is a finite morphism of degree 2
(see (6.6.7)(i)).

\bigskip

\centerline{\bf 1. The two Fano families ${\cal F}_{+}$ and }
\centerline{\bf ${\cal F}_{-}$ on $X \rightarrow {\bf P}^{2}.$}

\smallskip

{\bf (1.1)}
{\it The conic bundle surfaces $S_{l}$ and their
relatively minimal models $S(L)$.}

\smallskip

Let $p:X \rightarrow {\bf P}^{2}$ be a general
smooth conic bundle,
let
${\Delta} = \{ x \in {\bf P}^{2} : \mbox{ the fiber }
p^{-1}(x) \mbox{ is singular } \} \subset {\bf P}^{2}$
be
the discriminant curve of $p$, and let
${\pi} : \tilde{\Delta} \rightarrow {\Delta}$
be the double covering induced by the splitting
$p^{-1}(x) = {\bf P}^{1} \cup {\bf P}^{1}$ of the
fiber ${p}^{-1}(x), x \in {\Delta}$.
Let $d = deg({\Delta})$.
Since $X$ is assumed to be general, we can suppose that
$\Delta$ is a general plane curve of degree $d$
(see [Isk], Lemma 1 (iv)).
In addition, we shall always assume that
$d = deg({\Delta}) \ge 3$,
disregarding the trivial cases $d = 1,2$.

The general line $l \subset {\bf P}^{2}$ intersects
$\Delta$ in $d$ points $(x_{1},...,x_{d}).$
Let $S_{l} = p^{-1}(l).$  Then $S$ is a non-minimal ruled
surface over $l \cong {\bf P}^{1}$.  Indeed, $S$
inherits the conic bundle structure $p:S \rightarrow l$
from $p:X \rightarrow {\bf P}^{2}.$
Let $p^{-1}(x_{i}) = l_{i} + \overline{l}_{i}$,
$i = 1,...,d$ be the singular fibers of $S$. Since
$l$ is general, $S_{l}$ is smooth, and $l_{i}$ and
$\overline{l}_{i}$ are (-1)-curves on $S_{l}$.
The effective divisor
$L = l_{1} +...+ l_{d} \in S^{d}({\Delta})$
defines a morphism
${\sigma}(L): S_{l} \rightarrow S(L)$,
${\sigma}(L)$ =
${\bf{\Pi}}_{i = 1,...,d} \  {\sigma}(\overline{l}_{i})$,
where ${\sigma}(\overline{l}_{i})$ is the blow-down
of $\overline{l}_{i}$. The map $p$ defines on $S(L)$ a
structure of a ruled surface
$p(L): S(L) \rightarrow l$.

{\sc Definition.}
Let $p: S \rightarrow C_{0}$ be a ruled surface over the
curve $C$, and let $e = e(S)$ be the
{\it invariant} of the ruled surface $S$, i.e.
$e = min \{ C^{2}: C \mbox{ -- a section of } p \}$.

Call the section $C$ of $S$ -- a
{\it minimal section} of $S$ (resp. -- of $p$)
if $C^{2} = e(L)$.

By definition, the set of minimal sections of $S$ is
non-empty -- finite or not.

{\sc Definition.}
Call the minimal section $C$ -- an
{\it isolated minimal section}
if the set of minimal sections of $S$ is
finite.  Call the minimal section $C$ -- a
{\it isolated minimal section}
if the set of minimal sections of $S$ is
infinite.

\smallskip

Let
$\pi: \tilde{\Delta} \rightarrow \Delta$
be double
covering induced by $p$, and let

${\pi}_{d} = S^{d}({\pi}) :
S^{d}(\tilde{\Delta}) \rightarrow S^{d}({\Delta})$
be the $d^{th}$ symmetric power of $\pi$.
We identify $\mid {{\cal O}_{{\bf P}^{2}}}(1) \mid$
and $\mid {{\cal O}_{\Delta}}(1) \mid$ by the natural
isomorphism
$\cap : {{\bf P}^{2}}^{*} \rightarrow S^{d}({\Delta})$,
${\cap}(l) = l.{\Delta}$.

\bigskip

{\bf (1.2)}
{\it The components ${F}_{+}$ and ${F}_{-}$ in
$S^{d}(\tilde{\Delta})$, and their subsets
${F}^{0}_{e}$.}

\smallskip

Let
$F = {{\pi}_{d}}^{-1}(\cap {{\bf P}^{2}}^{*})
\subset S^{d}(\tilde{\Delta})$. It is well-known
(see e.g. [B3]) that the 2-dimensional
family $F$ splits into two connected components:
$F = F_{+} \cup F_{-}$.
We shall find two natural families ${\cal F}_{+}$ and
${\cal F}_{-}$ of curves on $X$
which correspond to the components $F_{+}$ and $F_{-}$.

Let $L \in F$, and let $l \in {{\bf P}^{2}}^{*}$ be the unique
line such that ${\cap}(l) = {\pi}_{d}(L)$.
Let $U \subset {{\bf P}^{2}}^{*}$ be the open subset:

$U = \{ l : l$ intersects ${\Delta}$
in a set of $d$ disjoint points,
and $S_{l} = p^{-1}(l)$ is smooth $\}$.

Let $F^{0} \subset F$ be the open subset
$F^{0} = {{\pi}_{d}}^{-1}(U)$. The map
$e: L \rightarrow e(L) = e(S(L))$, (see (1.1)), is well
defined on $F^{0}$. We introduce the subsets

$F^{0}_{e} = \{ L \in F^{0} : e(L) = e \} \subset F^{0}$.

We shall prove the following

\bigskip

{\bf (1.3)}
{\bf Lemma.}
{\sl
There exists an open subset $V \subset U$ such that
${{\pi}_{d}}^{-1}(V) \subset F^{0}_{-1} \cup F^{0}_{0}$.
}

{\bf Proof.}
Let $e^{+} = max \{ e(L): L \in F_{0} \}$.
Clearly,
$e^{+} \le 0$, and we have to see that $e^{+} = 0$.

Assume that $e^{+} \le -1$, and let $L \in F^{0}_{e}$.
Then the rational ruled surface $S(L) \cong {\bf F}_{-e}$,
and $S(L)$ has only one minimal section -- a section
$C$ such that $(C^{2})_{S(L)} = e$.
Let $L = l_{1} + ... + l_{d}$, and let
$L_{i} = L + \overline{l}_{i} - l_{i}, i \in \{ 1,...,d \}$.
Since $e(L) = e^{+}$ is maximal, $e_{i} = e(L_{i}) \le e^{+}$.
The ruled surface $S(L_{i})$ is obtained from $S(L)$ by
an elementary transformation $elm_{i}$ -- centered in
the point ${\sigma}(\overline{l}_{i}) \in S(L)$. In particular,
$\mid e_{+} - e_{i} \mid = 1$. Since $e_{i} \le e_{+}$,
$e_{i} = e^{+} - 1$. Therefore the unique minimal section
$C = C(L) \subset S(L)$ passes through the point
${\sigma}(\overline{l}_{i})$ (see e.g. [LN, Lemma (4.3)]).
Denote by $C, \  C \subset X$
also the proper preimage of $C \subset S(L)$ in
$S_{l} \subset X$. It follows that $C \subset X$
intersects the component $\overline{l}_{i}$,
and we have to see that $C$ does not intersect
the component $l_{i}$. In fact, if $C$ intersects also
the component $l_{i}$, then the line $l$ must be tangent
to $\Delta$ in the point
$x_{i} = p(l_{i}) = p(\overline{l}_{i})$, which is
impossible since $l \in U$.
If we repeat the last argument for any
$i \in \{ 1,...,d \}$, we obtain that (the proper preimage
$C \subset X$ of) the unique minimal section $C \subset S(L)$
intersects any of the components $\overline{l}_{i}$,
$i = 1,...,d$, and $C$ does not intersect $l_{i}$ for
$i = 1,...,d$. Clearly, $C \subset X$ is also the proper
preimage of the unique minimal section $C_{i}$ of
$S(L_{i})$, $i = 1,...,d$. Let
$\overline{L} = \overline{l}_{1} + ... + \overline{l}_{d}$.
It follows from the preceding that the map
$\overline{\sigma} : S_{l} \rightarrow S(\overline{L})$
sends $C$ isomorphically onto the unique minimal section
of $S(\overline{L})$, and
$e(\overline{L}) = e^{+} - d$.

Note that the same arguments as above imply that the subset
$F^{0}_{{e}^{+}}$ must contain an open subset of $F^{0}$.
In fact, the image of the map
$e: F^{0} \rightarrow {\bf Z}_{\le 0}$
is discrete, and $e(L)$ must be constant on an open subset
of any of the components of $F^{0}$. Let $e_{0}$ be the minimal
among
these values of $e$ for which $e(L) = e_{0}$ on an open set.
Just as above, if
$L = l_{01} + ... + l_{0d}$ $\in F^{0}_{e_{0}}$
is general, then the the minimal section $C_{0}$ of
$S(L_{0})$ does not intersect any of the components
$\overline{l}_{0i}$, $i = 1,...,d$. In particular,
$e(L_{0} + \overline{l}_{0i} - l_{0i}) = e_{0} + 1$,
for any $i = 1,...,d$. Therefore  $e = e_{0} + 1$
over an open subset of $F^{0}$. By repeating the same
argument we obtain by induction that
$e = e^{+}$ over an open subset $V$ of $F^{0}$.

Now, the assumption $e_{0} < 0$ implies that $e(L) = e^{+}$
for exactly one
$L \in {{{\pi}_{d}}}^{-1}(l)$, $L \in V$.
Indeed,
if $L = l_{1} + ... + l_{d} \in V$
and $L \in F^{0}_{{e}^{+}}$,  then the combinatorial
arguments from above imply:
$e(L^{'}) = e^{+} - d + \# \{ L \cap L^{'} \} < e^{+}$, for
any
$L^{'} \in {{\pi}_{d}}^{-1}({\pi}_{d}(L)) - \{ L \})$.
Therefore  the map
${\pi}_{d}: F \rightarrow {{\bf P}^{2}}^{*}$
admits an inverse, over the open subset
$V \subset F^{0}$ (by construction, $V$ is a subset
of $F^{0}_{{e}^{+}}$).

It follows that the closure of $V$, in $F^{0}_{{e}^{+}}$,
contains a 2-dimensional rational family.
Therefore the map
$L \rightarrow l$, where $l$ is the unique line such that
${\pi}_{d}(L) = {\Delta}.l$,  defines a linear system
of degree $d$ and of dimension $2$ on $\tilde{\Delta}.$
The uniqueness of the inverse $L \in V$ of $l$ implies
that $\tilde{\Delta}$ must be  projected
onto a plane curve of degree $d$, and the degree of the
projection is 1. In particular,
$g( \tilde{\Delta}) \le g({\Delta})$ -- contradiction.
Therefore, $e^{+}$ cannot be negative, i.e. $e^{+} = 0$.

Let $F^{0}_{+} = \{ L \in F : e(L) = 0 \}$. Clearly,
$F^{0}$ contains an open subset of $F$;  and we can
assume, without any confusion, that $F^{0}$ is open.
Since $S(L) \cong {\bf F}_{0}$ for $L \in F^{0}_{+}$,
the ruled surface $S(L) \rightarrow l$ has an infinite
set of minimal sections. In fact, $S(L)$ is a quadric.
Let $L_{i}$ be as above. Then $e(L_{i}) = -1$.
Therefore, there exists an open subset
$F^{0}_{-} \subset F$ such that $e(L) = -1$ for any
$L \in F^{0}_{-}$. {\bf q.e.d.}

\bigskip

{\bf (1.4)}
{\bf Corollary.}
{
Among the sets $F^{0}_{e}$,
$F^{0}_{0}$ and $F^{0}_{-1}$ are the only sets
in $F^{0}$ which contain open subsets of
$F = {\pi}_{d}^{-1}( \cap {{\bf P}^{2}}^{*} )$.
}

In fact, $F^{0}_{+} \cap F^{0}_{-}$ dominate over
an open subset of ${{\bf P}^{2}}^{*}$.
Therefore, the set
$U = \{ l \in {{\bf P}^{2}}^{*} : e(L) \in \{ 0,1 \}
\mbox{ for any } L \in {{\pi}_{d}}^{-1}(l.{\Delta}) \}$
contains an open subset of the plane.
{\bf q.e.d.}

\bigskip

{\bf (1.5)}
{\it
The Fano families of minimal sections
${\cal F}_{+}$ and ${\cal F}_{-}$.
}

\smallskip

Let $F_{+}$ and $F_{-}$ be the closures of $F^{0}_{+}$
and $F^{0}_{-}$ in $F$.
It follows from (1.2) and (1.3) that $F = F_{+} \cup F_{-}$,
i.e., outside a closed subset of codimension $\ge 1$,
the invariant $e(L)$ is either $0$ or $-1$.
If $L$ is general and $e(L) = -1$ then $S(L) \cong {\bf F}_{1}$,
and $S(L)$ has exactly one minimal section. Moreover
(see the proof of (1.2)), the proper preimage $C \subset X$ of
this minimal section does not intersect any of the components
$\overline{l}_{i}$ conjugate of the elements $l_{i}$ of $L$,
$i = 1,...,d$. Therefore, the proper preimage $C$ coincides
with the preimage of the minimal section (the preimage does
not contain components of degenerate fibers of
$p:X \rightarrow {\bf P}^{2}$). Therefore  we can define

${\cal F}_{-}$ = (the closure of)
$\{ C \subset X : C = C(L)$ is (the isomorphic preimage
of) the unique minimal section of $S(L)$  for some
$L \in F^{0}_{-} \}$.

Similarly,
we define
${\cal F}_{+}$ = (the closure of)
$\{ C \subset X : C \mbox{ is a minimal section of } S(L)
 \mbox{ for some } L \in F^{0}_{+} \}$.

Clearly, the family ${\cal F}_{-}$ is 2-dimensional, and
$p$ a defines the
the finite map
$p_{-} : {\cal F}_{-} \rightarrow {{\bf P}^{2}}^{*}$
of degree $deg(p_{-}) = 2^{d-1}$.
The family ${\cal F}_{+}$ is 3-dimensional, and $p$
defines a map
$p_{+} : {\cal F}_{+} \rightarrow {{\bf P}^{2}}^{*}$.
The general fiber ${p_{+}}^{-1}(l)$ consists of a finite
number (=$2^{d-1}$) of smooth rational 1-dimensional families
of curves on $X$.  Keeping in mind all this we postulate:

{\sc Definition.}
Call:
1. the 3-dimensional family ${\cal F}_{+}$ --
{\it the Fano family of non-isolated minimal sections of $p$}
;
2. the 2-dimensional family ${\cal F}_{-}$ --
{\it the Fano family of isolated minimal sections of $p$}.

\bigskip

\centerline{\bf  2. The versality of the family}
\centerline{\bf  of conic bundle surfaces  ${\cal S}[k]$.}

\smallskip

{\bf (2.1)}
{\it
Remarks.
}

Let $F_{+}$ and $F_{-}$ be the components of $F$, and let
$F_{+}^{reg} = F_{+} - Sing \ F_{+}$ and
$F_{-}^{reg} = F_{-} - Sing \ F_{-}$.
Consider e.g. the component $F_{-}$. Let
$Z \subset F_{-}$ be a smooth algebraic curve, and
assume that $Z$ is otherwise general; in particular,
we can assume that the general point $z \in Z$
belongs to $F_{-}^{reg}$.
The general point $z \in Z$ represents an effective
divisor $L = L(z) \in S^{d}(\tilde{\Delta})$ such that
${\pi}_{d}(L) = l$ is a line in ${\bf P}^{2}$.
Let
$p_{L} : S(L) \rightarrow l \cong {\bf P}^{1}$
be the ruled surface defined by $L$.
This way, the curve $Z$ defines the algebraic family of
ruled surfaces
${\cal S} \rightarrow Z$, such that the general fiber
$S(z) = S(L(z))$ is isomorphic to
${\bf F}_{1} =
{\bf P}({\cal O}_{{\bf P}^{1}} \oplus {\cal O}_{{\bf P}^{1}}(1))$.

The same arguments can be repeated for $Z \subset F_{+}$.
It follows from section 1 that the general fiber $S(z)$ of the
corresponding families
${\cal S} \rightarrow Z$ is a ruled surface of type
${\bf F}_{0}$ -- a smooth quadric with a fixed ruling.
Remember that the general member of a versal deformation
of a rational ruled
surface must be either of type ${\bf F}_{1}$ or of type
${\bf F}_{0}$. Moreover, two rational ruled surfaces
$S_{1}$ and $S_{2}$ can be deformed into each other
{\it iff} their invariants $e(S_{1})$ and $e(S_{2})$ have
the same parity (see [Se, Th.5, Th.13]). Now, the results
from section 1 imply that the general members of the families
${\cal S} \rightarrow Z$ are members of a versal deformation
of a rational ruled surface.

\bigskip

{\bf (2.2)}
{\it Versal deformations of 2-dimensional conic bundles.}

Let ${\rho} : {\cal S} \rightarrow Y$ be an algebraic family of
surfaces such that the general element $S(y) = {\rho}^{-1}(y)$
has a structure of a smooth conic bundle
$p(y) : S(y) \rightarrow C(y)$
over a (smooth) curve $C(y)$ of genus $g$, with $d$ degenerated
fibers.

Clearly, the general surface $S(y)$ admits exactly $2^{d}$
morphisms ${\sigma} : S(y) \rightarrow S_{\sigma}$,  where $S$
is a ruled surface
$p_{\sigma} : S_{\sigma} \rightarrow C(y)$,
such that
$p(y) = {\sigma} \circ p_{\sigma}$.
The map ${\sigma} \rightarrow e(S_{\sigma})$ (= the invariant
of the ruled surface $S_{\sigma}$) defines a map

$e : Y \rightarrow \{ \mbox{the subsets of } {\bf Z} \}$.

{\sc Definition.}
{\it Isolated and non-isolated minimal sections
of a conic bundle surface.}

Let $p:S \rightarrow C_{0}$ be a smooth conic
bundle surface over the curve $C_{0}$.
Let $C \subset S$ be a non-singular section of $p$
(i.e. $p$ maps $C$ isomorphically onto the base $C_{0}$,
and $C$ does not intersect the finite set of singular points
of the degenerate fibers of $p$),
and let $p(C):S_{C} \rightarrow C_{0}$ be the ruled surface
defined by $C \subset S$ (i.e. $S \rightarrow S(C)$ is the
contraction of the components
of the singular fibers of $p$ non-intersecting $C$, and
$p(C)$ is the induced ${\bf P}^{1}$-bundle map).

Call $C$ --
{\it a minimal section of the conic bundle surface $S$}
if $C$ is a minimal section of the ruled surface $S(C)$.
Call the minimal section $C$ of $S$ {\it isolated}
(resp. -- {\it non-isolated}) if $C$ is an isolated
(resp. -- non-isolated) minimal section of $S(C)$.

As it follows from [Se, Th. 13]
(see also (5.2) -- (5.7))
the general member of a versal
deformation of a ruled surface, over a curve of genus $g$,
is a ruled surface $S$ with invariant $e \in \{ g-1 , g \}$, and

(i). if $g = 1$ and $e = 0$ then $S$ is decomposable;

(ii) if $g = 1$ and $e = 1$, or if $g \ge 2$, then
$S$ is indecomposable.

In particular (see (5.2) -- (5.7)):
If $C$ is a minimal section of
the general versal element $S$,
then $C$ is isolated or not --
depending on $e = e(S) = g - 1$ or $e = g$.
This suggests the following

{\sc Definition.}
{\sl
Call the algebraic family of conic bundles
${\rho} : {\cal S} \rightarrow Y$
versally embedded, if the invariant $e$ of the general
$S_{\sigma}$ belongs to $\{ g-1 , g \}$
and  $S_{\sigma}$  fulfills (i), (ii).

Equivalently, the family ${\cal S}$ is versally embedded
if there exists an open subset
$U \subset Y$ such that
$e(U) \subset \{ g-1 , g \}$, and the ruled surfaces
$S_{\sigma} = {\sigma}(S(u)), u \in U$,  fulfill
(i),(ii).
}

In particular (see (5.2) -- (5.7)):
If $C$ is a (non-singular) minimal section
of the general versal element $S$,
then $C$ is isolated or not  --
depending on
$e = e(S(C)) = g - 1$ or $e = g$.

\bigskip

{\bf (2.3)}
{\bf Corollary.}
{\sl
Let $X$ be a smooth threefold,  which admits a regular map
$p : X \rightarrow {\bf P}^{2}$ which is a conic bundle
structure on $X$. Then the family of conic bundle surfaces
${\cal S} \rightarrow {{\bf P}^{2}}^{*}$, with fibers
$S(l) = p^{-1}(l)$, is versally embedded.
}

\bigskip

Let
$\mid {\cal O}_{{\bf P}^{2}}(k) \mid$
be the family of
plane curves of degree $k$, $k < d = deg({\Delta})$.
The intersection defines an isomorphism
$\cap : \mid {\cal O}_{{\bf P}^{2}}(k) \mid
\rightarrow \mid {\cal O}_{\Delta}(k) \mid$.

Let ${\cal C}[k] =
{{\pi}_{kd}}^{-1}( \mid {\cal O}_{{\bf P}^{2}}(k) \mid )$ =
${{\pi}_{kd}}^{-1} \circ
{\cap}^{-1}(\mid {{\cal O}_{{\bf P}^{2}}}(k) \mid$.
It follows from the definition that
${\cal C}[k]$ is a $(k+1).(k+2)/2 \  - 1$-dimensional subset
of $S^{k}(\tilde{\Delta})$, since the natural map
${\pi}_{d} : {\cal C}[k] \rightarrow
\mid {{\cal O}_{{\bf P}^{2}}}(k) \mid
\cong {\bf P}^{(k+1).(k+2)/2 \ -1}$
is finite of degree $2^{kd}$.

Let $C \subset {\bf P}^{2}$ be a general plane curve of
degree $k$, and let $S_{C} = p^{-1}(C)$. The conic bundle
structure $p : X \rightarrow {\bf P}^{2}$ defines a conic
bundle structure $p(C) : S(C) \rightarrow C$, on the
surface $S(C)$.

This way we obtain the family of conic
bundle surfaces
${\cal S}[k] \rightarrow
\mid {{\cal O}_{{\bf P}^{2}}}(k) \mid$,
parameterized by the space of plane curves of degree $k$.

\bigskip

{\bf (2.4)}
{\bf Theorem.}
{\sl
Let $p : X \rightarrow {\bf P}^{2}$ be a smooth conic bundle.
Let ${\pi} : \tilde{\Delta} \rightarrow {\Delta}$ be the
corresponding double covering of the discriminant curve
(${\Delta}$ is assumed to be smooth, and $\pi$ -- unbranched
and non-trivial). Let ${\cal F}_{-}$ be the 2-dimensional
family (the Fano surface of $X$) defined in section 1, and
assume that the following is true:

{\bf (*).}
The curves of the Fano surface ${\cal F}_{-}$ sweep $X$ out.

Then the family
${\cal S}[k] \rightarrow \mid {\cal O}_{{\bf P}^{2}}(k) \mid$
is versally embedded, for any positive integer
$k < d = deg({\Delta})$.
}

\bigskip

\centerline{\bf 3. Proof of Theorem (2.4).}

\smallskip

{\bf (3.1)}
{\sl
The level families ${\cal F}_{-}(m)$ and ${\cal F}_{+}(m)$.
}

Let $S_{\Delta} = p^{-1}({\Delta})$. Obviously,
the map $p : X \rightarrow {\bf P}^{2}$
defines an isomorphism

$p_{\Delta} : Sing(S_{\Delta}) \rightarrow
{\Delta} \subset {\bf P}^{2}$.

The inverse $s$ of $p_{\Delta}$ is  \cal{the Steiner map}:
$s : {\Delta} \rightarrow Sing(S_{\Delta})$,
$s : x \mapsto sing \ {p^{-1}(x)}$.

Let ${\cal F}_{+}$ and ${\cal F}_{-}$ be the
Fano families defined in section 1. The general
element  $C \in {\cal F}_{-}$  does not intersect
$s({\Delta})$ and defines by intersection
with $S_{\Delta} = p^{-1}({\Delta})$ the effective
divisor without multiple points
$L = L(C) \in F_{-}$.
Let $l \subset {\bf P}^{2}$ be the
line defined by $L$.
The morphism
${\sigma}(L) : S_{l} = p^{-1}(L) \rightarrow S(L)$
defined by $L$ (see section 1) maps $C$ isomorphically
onto the (-1) section $C(L)$ of $S(L) \cong {\bf F}_{1}$.
Since $C$ does not intersect $s({\Delta})$,
${(C^{2})_{S_{l}}} = {(C(L)^{2})_{S(L)}} = -1$, i.e.
$C$ is a (-1)-section also on the conic bundle surface
$S_{l}$.

Similar arguments imply that the general
$C \in {\cal F}_{+}$ is a 0-section of the conic bundle
surface $S_{l}$ defined by $L$.

Let $S_{l} = p^{-1}(l)$ be as above, and let
$C \subset S_{l}$ be a curve such that
$p : C \rightarrow l$ is an isomorphism. We call
such a curve a {\it section of the conic bundle
surface} $S_{l}$. Call $C$ {\it nonsingular} if
$C \cap s({\Delta}) = \oslash$.
Any nonsingular section $C$ of $S_{l}$ defines --
via an intersection with $S_{\Delta}$ -- a nonsingular
element $L$ of ${\cal F}$, i.e. $L$ has no multiple
points. The curve $C \subset X$ can be regarded also
as a section of the ruled surface $S(L)$.
Let $e(L) = e(S(L))$ be the invariant of $S(L)$.
Since $C$ is a nonsingular section of $S_{l}$,
the intersection number ${(C^{2})}_{S_{l}}$ is well
defined, and
${(C^{2})}_{S_{l}}$ $\equiv 0 ( \ mod \ e(L))$. Since
$e(L) \in \{ 0 , -1 \}$, the parity of
$C^{2}$ defines uniquely $e(L)$.
Let $C_{0}$ be a minimal section of $S(L)$. The general
choice of $l$ implies that $C_{0}$ can be considered as a
nonsingular section of $S_{l}$.

Denote by $f$ the general fiber of $S(L)$. The fiber $f$
can be considered as a fiber ( a conic)
$q = p^{-1}(x), x \in l - (l \cap {\Delta})$.
In particular, $C$ can be considered both: as an element of
the linear system $\mid C(L) + m.f \mid$ on $S(L)$,
and as an element of the linear system
$\mid C + m.q \mid$ on $S_{l}$; here $m = m(L)$ is defined by
$2m = {(C^{2})}_{S_{l}}$ - $e(L)$.
This suggests
to define the families

${\cal F}_{-}(m)$ =  the closure of
$\{ C$ -- a curve on $X:$
$C$ is a nonsingular section of some (nonsingular)
$S_{l}$, and ${(C^{2})}_{S_{l}}$ = $2m - 1 \}$.

We call the elements $C$ of
${\cal F}_{-}(m)$ -- {\it odd}
sections (of level $m$).

The same definition, up to parity, is used to define
the families ${\cal F}_{+}(m)$
of {\it even} Fano sections of level $m$
(the number $C^{2}$, for a general
even section $C$ of level $m$, is $2m$).

Let $l$ be a line in ${\bf P}^{2}$. It follows from the
preceding that:

\smallskip

{\sl
Any section $C$ of $S_{l}$
(i.e., a curve on $S_{l}$
which is mapped isomorphically to $l$ by $p$)
is an element of some of the level families
${\cal F}_{-}(m)$ or ${\cal F}_{+}(m)$.
}

\bigskip

{\bf (3.2)}
{\bf Corollary.}
{\sl
The conic bundle projection
$p : X \rightarrow {\bf P}^{2}$
defines, for any nonnegative integer $m$,
the surjective morphisms
$p_{-}(m) : {\cal F}_{-}(m) \rightarrow {{\bf P}^{2}}^{*}$
and
$p_{+}(m) : {\cal F}_{+}(m) \rightarrow {{\bf P}^{2}}^{*}$.

Moreover, if $d = deg({\Delta})$ and
$l$ is a general line in ${\bf P}^{2}$ then the
fiber $(p_{-}(m))^{-1}(l)$ has $2^{d-1}$ irreducible
components, and any such a component is isomorphic
to the projective space ${\bf P}^{2m}$.
Similarly, $(p_{+}(m))^{-1}(l)$ has $2^{d-1}$
irreducible components isomorphic to ${\bf P}^{2m+1}$.
}

{\bf Proof.}
Consider the fiber of $p_{-}(m)$.
Let $C \in (p_{-}(m))^{-1}(l)$ be general (esp. nonsingular),
and let $L = L(C) \in F$ be the divisor defined by $C$.
The minimal section $C(L)$ is nonsingular.
Therefore  $C(L) \subset S(L)$
can be regarded also as a curve on $X$; in particular

$C \in \mid C(L) + m.q \mid \cong \mid C(L) + m.f \mid$
$\cong {\bf P}^{2m}$.

Clearly, the nonsingularity of $C$ implies the uniqueness
of the minimal section $C(L)$ in the representation
$C \in \mid C(L) + m.q \mid$. Therefore, the components
of $(p_{-}(m))^{-1}(l)$ are in (1:1)-correspondence with the
$2^{d-1}$ sections of ${\cal F}_{-}$ "above" $l$, i.e. --
with the $2^{d-1}$ elements of $(p_{-}(0))^{-1}(l)$.

Similarly -- for $(p_{+}(m))^{-1}(l)$; note that
$(p_{+}(0))^{-1}(L)$ consists of $2^{d-1}$ components,
and any of these components is isomorphic to a ruling
of a quadric defined by a nonsingular minimal section
$C(L)(t) \in {\cal F}_{+}$ over $l$.

\bigskip

{\bf (3.3)}
{\it
Minimal sections of degenerate conic bundle
surfaces in $X$.
}

We shall use the level families in the study of the
minimal sections of the singular conic bundle surfaces
$S_{C} = p^{-1}(C)$, where $C$ is a reduced plane curve
of degree $k$ all the component of which are lines.

Let ${\bf S}_{k} \subset
\mid {{\cal O}_{{\bf P}^{2}}}(m) \mid$
be the 2k-dimensional locus of plane curves of
degree $k$, all of the components of which are
(possibly multiple) lines.
We shall see, after an
appropriate completion of the definitions, that
the induced family of degenerate conic bundle surfaces
${\cal S} \rightarrow {\bf S}_{k}$ is versally embedded,
provided the condition (*) of (2.4) takes place.

Let $x$ be a point of ${\bf P}^{2}$, and let
$q(x) \subset X$ be the conic
$q(x) = p^{-1}(x)$. Let ${\sigma}(x)$ =
${\sigma}_{10}(x) = \{ l \in {{\bf P}^{2}}^{*} : x \in l \}$
be the "Schubert line" of lines through $x$. Let
${\cal F}_{-}(x) = {\cal F}_{-}(0)(x)$ be the set of all
the odd Fano curves $C$ of level $0$ such that the line
$l = p(C)$ passes through $x$. Clearly,
${\cal F}_{-}(x) =
\cup \{ {\cal F}_{-}(l) : l \in {\sigma}(x) \}$, where
${\cal F}_{-}(l) =
\{ C \in {\cal F}_{-} : p(C) = l \}$.

If $l$ is general, then the set
${\cal F}_{-}(l)$ consists of $2^{d-1}$ elements:
$C_{1}(l), ... , {C_{{2}^{d-1}}}(l)$. Any of the
curves $C_{i}(l)$ intersects the conic $q(x)$
in a point ${\xi}_{i}(l) , i = 1, ..., 2^{d-1}$.

It follows from (2.4)(*) that any point
${\xi} \in q(x)$ can be represented in this way,
i.e., there exists $l \in {\sigma}(p({\xi}))$
such that ${\xi} = C \cap q(p({\xi}))$ for some
$C \in {\cal F}_{-}(p({\xi}))$.  It is easy to see
that (2.4)(*) also implies that the set
$( {\cal F}_{-} )_{\xi} =
\{ C \in {\cal F}_{-} : {\xi} \in C \}$
is finite for the general ${\xi} \in X$.

Let $C_{0} = n_{1} + ... + n_{k}$ be a general element
of ${\bf S}_{k}$, and let
${S_{C_{0}}} \subset X$ be the singular conic bundle surface
${S_{C_{0}}} = p^{-1}(C_{0})$. The surface $S_{C_{0}}$ has $k$
irreducible components -- the rational conic bundle
surfaces
$S_{i} = p^{-1}(n_{i}) , i = 1,...,k$.
Since $C_{0}$ is general, all the $S_{i}$ can be assumed
to be smooth.

{\it Definition.}
Call a {\it section} of $S_{C_{0}}$ any {\it connected} curve
$C = C_{1} + ... + C_{k} \subset X$ such that
$C_{i}$ is a section of $S_{i}$.
Respectively, we call the section $C$ {\it nonsingular}
if $C_{i}$ is nonsingular, $i = 1,...,k$. Clearly, the
nonsingular section $C$ can be singular as a curve.

Since $C_{0}$ is general, we can assume that any of the
ruled surfaces $S_{i}({\sigma})$ is either
${\bf F}_{0}$ or ${\bf F}_{1}$.
We shall find the minimal sections of $S_{C_{0}}$.

Fix an order, say $[1,2,...,k]$. Let
$C_{1} \in {\cal F}_{-}$ be any of the $2^{d-1}$
odd Fano sections of level $0$ such that
$p(C) = n_{1}$. Let $x_{12} = n_{1} \cap n_{2}$, and
let $q(x_{12}) = p^{-1}(x_{12})$. The general choice
of $C_{0}$ implies that $q(x_{12})$ is smooth, and
that (because of (2.4)(*))
${\xi}_{12} = C_{1} \cap q(x_{12})$ is a general point
of $q(x_{12})$. It follows from (3.2) that the set
$\{ C \in {\cal F}_{+} :
p(C) = n_{2}, {\xi}_{12} \in C \}$ is finite
(it has $2^{d-1}$ =
$ \# \{ F^{o} : F^{o} - \mbox{a component of }
{\cal F}_{+}(n_{2}) \}$
elements, any of which can be assumed to be a nonsingular
section of $S_{2}$).
Let $C_{2}$ be any of these curves. Obviously,
$C_{1} + C_{2}$ is a minimal section of
$S_{1} + S_{2}$, and the invariant of the corresponding
reducible ruled surface
$S(C_{1} + C_{2})$ = $S(C_{1}) + S(C_{2})$
is $(C_{1})^{2} + (C_{2})^{2} = -1 = p_{a}(n_{1} + n_{2}) - 1$;
i.e. --  the singular (reducible) ruled surface $S(C_{1} + C_{2})$
is of type ${\bf F}_{1}$.
We shall prove
by induction the following:

\smallskip

{\bf (3.3.1)}
{\bf Lemma.}
{
If $C_{0} = n_{1} + ... n_{k}$ is general, then
the surface $S_{C_{0}}$ has only a finite number of sections
$C_{1} + ... + C_{k}$ such that

(i). if $1 \le i \le k$ and $i \in 2{\bf Z}$ then
$C_{i} \in {\cal F}_{-}(i/2)$;

(ii). if $1 \le i \le k$ and $i \in 2{\bf Z} + 1$ then
$C_{i} \in {\cal F}_{+}((i-1)/2)$.
}

{\bf Proof.}
Assume that (3.3.1) is true for $i = 1,2,...,k-1$.
Let $x_{ik} = n_{k} \cap n_{i}$,
let $q(x_{ik}) = p^{-1}(x_{ik})$, and
let ${\xi}_{ik} = C_{i} \cap q(x_{ik})$,
$i = 1,2,...,k-1$.
Since
$C_{0}$ is general, we can assume that $n_{k}$
is a general line in ${\bf P}^{2}$, and
$(x_{1k},...,x_{k-1,k})$ is a general $(k-1)$-tuple
of points on $n_{k}$.
Moreover, since $n_{i} \in {\sigma}(x_{ik})$ are
general, the points
${\xi}_{ik} \in q(x_{ik})$ can be assumed to be general
(see above).

Let, for definiteness, $k$ be even. It follows from (3.2)
that there is only a finite number of sections
$C \in {\cal F}_{-}(k/2)$ such that $p(C) = n_{k}$,
which pass through the points
${\xi}_{1k},...,{\xi}_{k-1,k}$. Indeed, the general
$(k-1)$-tuple $({\xi}_{1k},...,{\xi}_{k-1,k})$
imposes $k-1$ independent conditions on the elements
of any of the the finite number of linear systems
$\mid C_{0} + (k/2).q \mid$ on $S_{k} = p^{-1}(n_{k})$,
where $C_{0}$ is a (-1)-section of $S_{k}$.

By induction, $e(C_{1} + ... + C_{k})$ =
$e(C_{1} +...+ C_{k-1}) + (C_{k})^{2}$ =
$(p_{a}(n_{1} + ... n_{k-1}) - 1) + (-1 + 2.k/2)$ =
$p_{a}(n_{1} + ...+ n_{k}) - 1$.

If $k$ is odd, the proof is similar.
Lemma (3.3.1)
is proved.

\smallskip
{\bf (*).}
{\it Definition} of {\bf level} and {\it weight}.

Let $C = C_{1} + ... C_{k}$ be a section of
$S_{1} + ... + S_{k}$, and let
$m_{i}$ be the level of $C_{i}$,
$i = 1,...,k$. Call the vector
${\bf level}(C) = (m_{1}, ... , m_{k})$
the {\it level} of $C$. Call the {\it weight}
of $C$ the number
$weight(l) = m_{1} + ... + m_{k}$.

In particular, if $C$ is a section
of type described in (3.3.1) and $k = 2k_{0}$ is even,
then $weight(C) = (k_{0} - 1).k_{0}$; if $k = 2k_{0}+1$ is
odd then $weight(C) = k_{0}^{2}$.

Call a section $C$ of $S_{C_{0}}$ {\it isolated} if
the linear system $\mid C \mid$ on $S_{C_{0}}$ is trivial.

\smallskip

{\bf (3.3.2)}
{\bf Corollary.}
{
If $C_{0} = n_{1} + .. n_{k} \in {\bf S}_{k}$ is general
then the surface $S_{C_{0}}$ has only a finite number of
isolated sections. These are the sections $C$ of $S_{C_{0}}$
for which

$weight(C) = (k_{0} - 1)k_{0}$ if $k = 2k_{0}$ is even;

$weight(C) = k_{0}^{2}$ if $k = 2k_{0} + 1$ is odd.
}

The proof of (3.3.2) is purely combinatorial, and similar
to the proof of (3.3.1), where the isolated sections of
$S_{C_{0}}$ for which $m_{2i} = m_{2i+1} = i, 1 \le i \le k$,
are described.

Let $C = C_{1} + ... C_{k}$ be a nonsingular section of
$S_{C_{0}}$. Call $C$ minimal if $C$ is a section of a minimal
self-intersection $C^{2} = C_{1}^{2} + ...C_{k}^{2}$ on
$S_{C_{0}}$.

It is not hard to see also if $C$ is an isolated section
of such a $S_{C_{0}}$ then $C$ is minimal and
$e(C) = p_{a}(n_{1} + ... + n_{k}) - 1$ =
$(k-1).(k-2)/2 -1$ = $e(S(C_{1} + ... + C_{k}))$.

By using just the same arguments as in the proof
of (3.3.1) (based on the assumption (2.4)(*)), we obtain
the following

\smallskip

{\bf (3.3.3)}
{\bf Corollary.}
{\sl
Let $C_{0} \in {\bf S}_{k}$ be general, and let $C$
be a non-isolated minimal section of $S_{C_{0}}$.
Then

$weight(C) = (k_{0} - 1).k_{0} + 1$, if $k = 2k_{0}$ is even;

$weight(C) = k_{0}^{2} + 1$, if $k = 2k_{0} + 1$ is odd.

Moreover, the set of non-isolated (continual) minimal sections
of (the general) $S_{C_{0}}$ form a 1-dimensional algebraic
family, the components of which are finite number of rational
curves.
}

\bigskip

{\bf (3.4).}
{\it End of the proof of Theorem(2.4).}

\smallskip

Assume that the family
${\cal S}[k] \rightarrow \mid {{\cal O}_{{\bf P}^{2}}}(k) \mid$
is not versally embedded.  By construction,
${\bf S}_{k} \subset \mid {{\cal O}_{{\bf P}^{2}}}(k) \mid$.
Let ${\cal S}_{k} \rightarrow {\bf S}_{k}$ be the restriction
of ${\cal S}_{k}$ on the preimage of ${\bf S}_{k}$.

Let $g = g(k) = (k-1).(k-2)/2 $ be the genus of the
general plane curve of degree $k$.  Then (3.3.1) - (3.3.3) imply:
$e(S) = \{ g-1 , g \}$, for the general
$S \in {\cal S}_{k}$.  Fix the general k-tuple of lines
$C_{0} = n_{1} + ... + n_{k}$ such that
$e(S(C_{0})) = \{ g(k) -1 , g(k) \}$.
Let $U \subset \mid {{\cal O}_{{\bf P}^{2}}}(k) \mid$ be
a sufficiently small analytic neighborhood of $C_{0}$.
The general curve $C \in U$ is a smooth plane curve of
degree $k$, and the invariant function
$e : U \rightarrow \{ \mbox{ the subsets of } {\bf Z} \}$
takes discrete values, on the elements  $C$ of the continuous
set $U$.  Now, standard arguments including the general
position of $C \in U$ imply that
$e(C) = e(S(C)) = e(C_{0}) = \{ g-1, g \}$
for the general element $C \in U$.
Theorem (2.4)
is
proved.

\bigskip

\centerline{\bf 4. The intermediate jacobian
$(J(X),{\Theta})$ as a Prym variety
$P(\tilde{\Delta} , {\Delta})$,}
\centerline{\bf and the family ${\cal S}[d-3]$.}

\bigskip

{\bf (4.1).}
{\it The Wirtinger description of $(J(X),{\Theta})$ as
a Prym variety, and the canonical families ${\cal C}_{+}$
and ${\cal C}_{-}$ of minimal sections.}

\smallskip

{\bf (4.1.1)}
Let $(\tilde{\Delta},{\Delta})$ be the discriminant pair of
$p:X \rightarrow {\bf P}^{2}$, and let
${\pi}: \tilde{\Delta} \rightarrow {\Delta}$ be the
induced double covering. Without any essential restriction
we assume that the discriminant curves are smooth, and
$\pi$ is unbranched.
It is well-known that the principally polarized intermediate
jacobian $(J(X),{\Theta})$ can be identified with the
Prym variety $P(\tilde{\Delta},{\Delta})$ defined by the double
covering ${\pi}: \tilde{\Delta} \rightarrow {\Delta}$
(see e.g.[B1]). Here we remember the Wirtinger description
of $P(\tilde{\Delta},{\Delta})$ by sheaves on $\tilde{\Delta}$
(see e.g. [W2]).

Let $d = deg({\Delta})$, and let $g = (d-1)(d-2)/2 = g({\Delta})$
be the genus of $\Delta$. The map $\pi$ induces the natural map
$Nm : {\bf Pic}(\tilde{\Delta}) \rightarrow
      {\bf Pic}({\Delta})$.

Let ${\omega}_{\Delta}$ be the canonical sheaf of $\Delta$.
Then the fiber ${Nm}^{-1}({\omega}_{\Delta})$ splits into
two components:

$P^{+} = \{ {\cal L} \in {\bf Pic}^{2g-2}(\tilde{\Delta}) :
            Nm({\cal L}) = {\omega}_{\Delta} ,
            h^{0}({\cal L}) \mbox{ even} \}$,

$P^{-} = \{ {\cal L} \in {\bf Pic}^{2g-2}(\tilde{\Delta}) :
            Nm({\cal L}) = {\omega}_{\Delta} ,
            h^{0}({\cal L}) \mbox{ odd}  \}$.

Both $P^{+}$ and $P_{-}$ are translates of the Prym variety
$P = P(\tilde{\Delta},{\Delta}) \subset
J(\tilde{\Delta}) = {\bf Pic}^{0}(\tilde{\Delta})$ ;
$P$ is the connected component of ${\cal O}$ in the
kernel of
$Nm^{0} : {\bf Pic}^{0}(\tilde{\Delta}) \rightarrow
          {\bf Pic}^{0}({\Delta})$

The general sheaf ${\cal L} \in P^{+}$ is
non effective, i.e. the linear system
$\mid {\cal L} \mid$ is empty. The set
${\Theta} = \{ {\cal L} \in P^{+}:$
$\mid{\cal L}\mid \ne \oslash \}$ =
$\{ {\cal L} \in P^{+} : h^{0}({\cal L}) \ge 2 \}$ is
a copy of the theta divisor of the p.p.a.v.
$P_{+} \cong P$. Since the general sheaf
${\cal L} \in P^{-}$ is effective, this suggests
to introduce the following two subsets
of $S^{2g-2}(\tilde{\Delta})$:

$Supp({\Theta}) =
\{ l \in \mid {\cal L} \mid : {\cal L} \in {\Theta} \}$,
   $Supp(P^{-}) =
\{ l \in \mid {\cal L} \mid : {\cal L} \in P^{-} \}$.

Clearly,
$dim \ Supp({\Theta}) = dim \ Supp(P^{-}) =
dim \ (P) = g - 1$. Indeed, the general fiber
${{\phi}_{\cal L}}^{-1} ({\cal L})$ of the
natural map
${\phi}_{\cal L} :
 Supp({\Theta}) \rightarrow {\Theta}$
coincides with the linear system
$\mid {\cal L} \mid \cong {\bf P}^{1}$,
and the general fiber of
${\phi}_{\cal L} :
 Supp(P^{-}) \rightarrow P^{-}$
is $\mid {\cal L} \mid \cong {\bf P}^{0}$.

The Gieseker-Petri inequalities for the
codimension of the special subsets of
$P^{+}$ and $P^{-}$ (see [W2])
imply that the general point
$L$ of $Supp({\Theta})$ lies on the
1-dimensional fiber over the general
point ${\cal L}$ of ${\Theta}$
(similarly -- for $Supp(P^{-})$).

We identify the effective sheaf ${\cal L}$
and the set of effective divisors
$L : L \in \mid {\cal L} \mid$.

Let
$S^{2g-2}(\pi) : S^{2g-2}(\tilde{\Delta})
     \rightarrow S^{2g-2}({\Delta})$
be the $(2g-2)$-th symmetric power of $\pi$, and let
$\mid {\omega}_{\Delta} \mid$
$\cong$  $\mid{\cal O}_{\Delta} (d-3) \mid$
$\cong \mid {\cal O}_{{\bf P}^{2}} (d-3) \mid$
$\cong {\bf P}^{g-1}$
be the canonical system of $\Delta$.
We shall use equivalently any of the different
interpretations of the elements of this system,
as it is written just above.

\bigskip

{\bf (4.1.2)}
{\it
The canonical families
${\cal C}_{+}$ and ${\cal C}_{-}$
of non-isolated and isolated minimal sections
of $p:X \rightarrow {\bf P}^{2}$.
}

Let ${\cal C}_{+}$ be the closure of the set

${\cal C}^{0}_{+}$ =
   $\{ C$
   -- a nonsingular section of
   $p: X \rightarrow {\bf P}^{2} :$
   $C_{0} = p(C)$
   is a curve from the canonical system for
   ${\Delta},$
   $\mbox{ } C_{0}$
   intersects transversely
   ${\Delta},$
   and $C$
   is a non-isolated minimal section of
   $S_{ C_{0} }  \}$.

Similarly, let
${\cal C}_{-}$ be the closure of
${\cal C}^{0}_{-} = \{ \mbox{ ... } ,
    \mbox{ and } C
    \mbox{ is an isolated (minimal) section of }
    S_{ C_{0} } \}$.

In the definitions of
${\cal C}_{+}$ and ${\cal C}_{-}$, we use
implicitly Theorem 2.4. In fact,
if the family $S[d-3]$ is not versally
embedded, then any minimal section
of the general conic bundle surface
$S_{ C_{0} }$ , $C_{0}$ -- a plane curve
of degree $d-3$,  will be isolated
(see e.g. [Se] or
[LN]).

Let $S_{\Delta} = p^{-1}({\Delta})$, let
${\psi}^{0} : {\cal C}_{+}^{0} \cup {\cal C}_{-}^{0}
 \rightarrow S^{2g-2}(\tilde{\Delta})$,
${\psi}(C) \mapsto L(C) = C \cap S_{\Delta}$, and let
${\psi} : {\cal C}_{+} \cup {\cal C}_{-}
\rightarrow S^{2g-2}(\tilde{\Delta})$
be the completion of ${\psi}^{0}$.

Denote by
$C_{+} = {\psi}({\cal C}_{+})$, and
$C_{-} = {\psi}({\cal C}_{-})$
the
${\psi}-$images of ${\cal C}_{+}$ and
${\cal C}_{-}$.

\bigskip

{\bf (4.2)}
{\bf Lemma.}
{\sl
The non-ordered pairs $\{ C_{+} , C_{-} \}$ and
$\{ Supp({\Theta}) , Supp(P^{-}) \}$
of subsets of $S^{2g-2}(\tilde{\Delta})$
coincide.
}

{\bf Proof.} It rests only to be seen that
the sets $C_{+} \cup C_{-} \subset
S_{2g-2}(\tilde{\Delta})$ and
${S^{2g-2}{\pi}}^{-1}(\mid {\omega}_{\Delta} \mid )$
coincide.
{\bf q.e.d.}

\bigskip

{\bf (4.3)}
{\it The Abel-Jacobi images of the families
${\cal C}_{+}$ and ${\cal C}_{-}$ }

\smallskip

Let
$J(X) = (H^{2,1}(X))^{*} / H_{3}(X,{\bf Z}) \ mod \mbox{ torsion}$

be the intermediate jacobian of $X$, provided with the principal
polarization ${\Theta}_{X}$ defined by the intersection of
3-chains on $X$. It is well known (see [B1]) that
$(J(X),{\Theta}_{X})$ is isomorphic, as a p.p.a.v., to the
Prym variety $(P,{\Theta})$ of the discriminant pair
$(\tilde{\Delta},{\Delta})$. Let

${\Phi}_{+} : {\cal C}_{+} \rightarrow J(X) \cong P$ and
${\Phi}_{-} : {\cal C}_{-} \rightarrow J(X) \cong P$

be the Abel-Jacobi maps for the families
${\cal C}_{+}$ and ${\cal C}_{-}$ of algebraically
equivalent 1-cycles on $X$.
Let
$Z_{+} = {\Phi}_{+}({\cal C}_{+})$ and
$Z_{-} = {\Phi}_{-}({\cal C}_{-})$
be the
images of ${\Phi}_{+}$ and ${\Phi}_{-}$.
We shall prove the following

\bigskip

{\bf (4.4)}
{\bf Theorem.}
{\sl
One of the following two alternatives always takes place:

{\bf (1).}
$h^{0}({\psi}(C)) = 2$  for the general
$C \in {\cal C}_{+}$  $\Leftrightarrow$
$h^{0}({\psi}(C)) = 1$  for the general
$C \in {\cal C}_{-}$,  and then

  (i).  $Z_{+}$ is a copy of the theta divisor ${\Theta}_{X}$,

  (ii). $Z_{-}$ coincides with $J(X)$;

{\bf (2).}
$h^{0}({\psi}(C)) = 1$  for the general
$C \in {\cal C}_{+}$  $\Leftrightarrow$
$h^{0}({\psi}(C)) = 2$  for the general
$C \in {\cal C}_{-}$,  and then

  (i).  $Z_{+}$ coincides with $J(X)$,

  (ii). $Z_{-}$ is a copy of the theta divisor ${\Theta}_{X}$.
}

{\bf Proof.}
The map
${\phi} = {\phi}_{\cal L} : Supp({\Theta}) \cup Supp(P^{-})
\rightarrow {\Theta} \cup P^{-}$ introduced above, can be
regarded as the (Prym)-Abel-Jacobi map for the sets
of algebraically equivalent (2g-2)-tuples of points
$Supp({\Theta}) \subset S^{2g-2}(\tilde{\Delta})$  and
$Supp(P^{-}) \subset S^{2g-2}(\tilde{\Delta})$,
to the Prym variety $P \cong J(X)$.

According to Lemma (4.2),
$C_{+} = {\psi}({\cal C}_{+})$ coincides either
with $Supp({\Theta})$,  or with  $Supp(P^{-})$.
Alternatively,
$C_{-} = {\psi}({\cal C}_{-})$ coincides either
with $Supp(P^{-})$, or with
$Supp({\Theta})$.

Let e.g. $C_{+} = Supp({\Theta})$. Then
$h^{0}({\psi}(C)) = 2$  for the general
$C \in {\cal C}_{+}$,
$h^{0}({\psi}(C)) = 1$  for the general
$C \in {\cal C}_{-}$; and we have to see
that $Z_{+} \cong {\Theta}$, and
$Z_{-} = J(X) \cong P$.

Let $C \in {\cal C}_{+}$ be general, and let
$z = {\Phi}_{+}(C) \in J(X)$  be the
Abel-Jacobi image of $C$.  Since $C$ is general,
$C$ is a nonsingular section of the conic
bundle surface $S_{p(C)} \subset X$, and the
effective divisor
$L = L(C) = {\psi}(C) \in Supp({\Theta})$ is well defined.

We can also assume that $p(C)$ is nonsingular, and
$p(C)$ intersects $\Delta$ transversely. In particular,
the effective divisor $L = L(C)$ does not contain
multiple points. We shall prove the following

\bigskip

{\bf (*).}
{\bf Lemma.}
{\sl
Let $C'$ and $C''$  $\in {\cal C}_{+}$ be such that
${\psi}(C') = {\psi}(C'') = L$, and let
$z' = {\Phi}_{+}(C')$, $z" = {\Phi}_{+}(C'')$.
Then $z' = z''$.
}

{\bf Proof of (*).}
Since ${\psi}(C') = {\psi}(C'')$, the curves $C'$
and $C''$ have the same p-image
$C_{0} = p(C') = p(C'')$, and $C'$ and $C''$ are
non-isolated sections of the conic bundle
surface $S_{C_{0}} = p^{-1}(C_{0})$.
Let $L = l_{1} + ... + l_{2g-2}$, and
$x_{i} = p(l_{i}), i = 1, ... , 2g-2$.
The degenerate fibers of
$p : S_{C_{0}} \rightarrow C_{0}$
are the singular conics
$q(x_{i}) = p^{-1}(x_{i}) = l_{i} + \overline{l_{i}}$.
By assumption $C'$ and $C''$  intersect simply
any of the components $l_{i}$, and does not intersect
any of $\overline{l_{i}}$.

Let $C$ be any nonsingular section of $S_{C_{0}}$ such that
${\psi}(C) = C \cap S_{\Delta} = L$, e.g. $C = C'$.  Then
${\bf Div}(S_{C_{0}}) =
p^{*}({\bf Div}(C_{0})) + {\bf Z}.l_{1} + ... + {\bf Z}.l_{2g-2} +
 {\bf Z}.C$.

Since $(C' - C'').q = 1 - 1 = 0$,
and $(C' - C'').l_{i} = 0, (i = 1,...,2g-2)$,
the divisor $C' - C''$
belongs to $p^{*}({\bf Div}(C_{0}))$; i.e.
$C' - C'' = p^{*}{\delta}$ for some
${\delta} \in {\bf Div}(C_{0})$.

Obviously,
$deg({\delta}) = 0$.
Represent $\delta$ as a difference of two
effective divisors (of the same degree):
${\delta} = {\delta}_{1} - {\delta}_{2}$.
Without loss of the generality we can assume that
the sets
$Supp({\delta}_{1})$  and  $Supp({\delta}_{2})$
are disjoint.
Therefore,
$p^{*}(C'-C'') =
p^{-1}({\delta}_{1}) - p^{-1}({\delta}_{2})$ is a
sum of fibers of $p$, with positive and negative
coefficients, and of total degree $0$.

Since all the fibers of $p : X \rightarrow {\bf P}^{2}$
are rationally equivalent, the rational cycle class
$[p^{-1}({\delta})]$,
of $p^{-1}({\delta})$, is $0$, in the Chow ring
$A.(X)$.  Since the Abel-Jacobi map for any family
of algebraically equivalent 1-cycles on $X$ factors
through the cycle class map, the curves $C'$ and $C''$
have the same Abel-Jacobi image, i.e. $z' = z''$.
This proves (*).

\smallskip

It follows from (*) that the Abel-Jacobi map ${\Phi}_{+}$
factors through $\psi$, i.e., there exists a well-defined
map
$\overline{\Phi}_{+} : Supp({\Theta}) \rightarrow Z_{+}$,
such that ${\Phi} = \overline{\Phi}_{+} \circ {\psi}$.

\smallskip

Let $C \in {\cal C}_{+}$ be general, and let
$L = L(C) = {\psi}(C)$.  Let
${\cal L} = {\phi}(L)$  be the sheaf defined
by the 1-dimensional linear system of effective divisors
linearly equivalent to $L$.
Let
${\cal C}_{+}({\cal L}) =
{\psi}^{-1} ( \mid {\cal L} \mid )$
be the preimage of  $\mid {\cal L} \mid$  in
${\cal C}_{+}$.
Since  ${\Phi}_{+}$  factors through  $\psi$, and
${\Phi}_{+}$ is a map to an abelian variety (the
intermediate jacobian  $J(X)$  of  $X$), the map
$\overline{\Phi}_{+}$ contracts rational subsets
of  $Supp({\Theta})$  to points.
However,
${\psi}({\cal C}_{+}({\cal L})) \cong
\mid {\cal L} \mid  \cong  {\bf P}^{1}.$
Therefore, there exists a point
$z = z({\cal L}) \in Z_{+}$  such that
${\Phi}_{+} ({\phi}^{-1}({\cal L}))$ =
${\Phi}_{+}({\cal C}_{+}({\cal L}))$ =
$\overline{\Phi}_{+}(\mid {\cal L} \mid)$ =
$ \{ z \} \subset Z_{+}$.

Clearly  $z = {\Phi}_{+}(C)$, and the
uniqueness of the sheaf ${\cal L}$  defined
by  $C$, implies that the correspondence
${\Sigma}$ =
$\{ (z , {\cal L}) :  z = {\Phi}_{+}(C),
{\cal L} = {\phi} \circ {\psi} (C) , C \in {\cal C}_{+} \}$
is generically (1:1).

Let  $i : {\Sigma} \rightarrow Z_{+}$  and
$j : {\Sigma} \rightarrow {\Theta}$
be the natural projections.  The general choice of
$C \in {\cal C}_{+}$, and the identity
${\psi}({\cal C}_{+}) = Supp({\Theta})$,
imply that  $j$  is surjective.  Therefore
$Z_{+}$  and  ${\Theta}$  are birational.
In particular, $Z_{+}$ is a divisor in
$J(X) \cong P$.  It is not hard to see that the map
$i \circ j^{-1} : {\Theta} \rightarrow Z_{+}$
is regular.  In fact,  let
${\cal L}$  be  any
sheaf which belongs to  $\Theta$.
The definition of  $\phi$  implies that
${\phi}^{-1}({\cal L})$ coincides with the
linear system  $\mid {\cal L} \mid$, which
is an (odd dimensional) projective
space. Therefore,  $\overline{\Phi}_{+}$  contracts
the connected rational set
${\psi}^{-1}({\cal L})$ to a unique point
$z = z(L)$, i.e. $i \circ j^{-1}$ is regular
in  ${\cal L}$. It follows that $Z_{+}$ is
biregular to the divisor of principal polarization
$\Theta$, i.e.  $Z_{+}$ is a translate of  $\Theta$.

The coincidence  $Z_{-} = J(X)$  follows in a similar way.

In case (2), the only difference is that the general fiber
of  ${\psi}$  is discrete, since the minimal sections
$C \in {\cal C}_{-}$  which majorate the general
$L \in Supp({\Theta})$, are isolated.
Theorem 4.4 is proved.

\bigskip

\centerline{\bf 5.  The fibers of the Abel-Jacobi maps
${\Phi}_{+}$ and ${\Phi}_{-}$.}

\smallskip

We shall describe explicitly the general fibers of
${\Phi}_{+}$ and ${\Phi}_{-}$ in either of the cases
(4.4)(1) and (4.4)(2).

\bigskip

{\bf (5.1)}
It follows from Theorem (4.4) that the fibers of ${\Phi}_{+}$
and ${\Phi}_{-}$ depend closely on the alternative conclusions:
$Z_{+} = {\Theta}$, or $Z_{-} = {\Theta}$.  The examples show
that any of the two alternatives (4.4)(1) -- (4.4)(2) can be
true, depending on the choice of the conic bundle
$p : X \rightarrow {\bf P}^{2}$ (see section 6).
Notwithstanding, theorems (2.4) and (4.4),  and some results
regarding versal deformations of ruled surfaces (see [Se])
and subbundles of rank 2 vector bundles over curves
(see e.g. [LN]),  make it possible to describe these fibers
in each of the cases (1) and (2).

\bigskip

{\bf(5.2)}
{\it Minimal sections of ruled surfaces and maximal subbundles
of rank 2 vector bundles on curves.}

\smallskip

Here we collect some known facts about ruled surfaces
and rank 2 vector bundles over curves.

\smallskip

Any ruled surface $S$ over a smooth curve $C$ can be
represented as a projectivization ${\bf P}_{C}(E)$
of a rank 2 vector bundle $E$ over $C$.  Clearly,
${\bf P}_{C}(E)$  is a ruled surface for any such $E$,
and ${\bf P}_{C}(E) \cong {\bf P}_{C}(E')$  {\it iff}
$E = E' \otimes {\cal L}$ for some invertible
sheaf ${\cal L}$;  here we identify vector bundles
and the associated free sheaves.

Call the bundle $E$ {\it normalized} if
$h^{0}(E) \ge 1$, but $h^{0}(E \otimes {\cal L}) = 0$
for any invertible ${\cal L}$ such that
$deg({\cal L}) < 0$ ( see [H, ch.5,\#2]).

The question is:

{\bf (*).} {\it How many normalized rank 2 bundles represent
the same ruled surface ?}

The answer depends on the choice of the curve $C$  (esp. --
on the genus $g = g(C)$ of $C$), and on the choice of the
ruled surface $S$ over $C$.
Let $p : S \rightarrow C$ be the natural fiber structure
on $S$. We shall reformulate the question (*) in the terms of
sections of $p$.

{\sc Definition.}
Call the section $C \subset S$ {\it minimal} if $C$ is
a section on $S$ for which the number $(C.C)_{S}$ is
minimal. Let $C$ be a minimal section of $S$.
The number $e = e(S) = (C.C)_{S}$ is an integer invariant
of the ruled surface $S$.  The number $e(S)$ coincides with
$deg(E) := deg(det(E))$, where  $E$ is any normalized
rank 2 bundle which represents $S$ ( i.e. -- such that
$S \cong {\bf P}_{C}(E)$) (see e.g. [H, ch.5,\#2]).
We call the number $e = e(S)$ the {\it invariant} of $S$.

{\bf Remark.}
Here, in contrast with the definition in use,
we let
$e(S) := - \mbox{ (the invariant of } S)$.

The new question is:

{\bf (**).} {\it How many minimal sections lie on the
same ruled surface ?}

The two questions are equivalent in the following sense:
Let $E$ be normalized and such that ${\bf P}(E) = S$.
By assumption  $h^{0}(E) \ge 1$. Therefore $E$ has
at least one section $s \in H^{0}(E)$. The bundle
section $s$ defines (and is defined by)
an
embedding
$0 \rightarrow {\cal O}_{C} \rightarrow E$.
The sheaf ${\cal L}$, defined by the cokernel of this
injection, is invertible, and ${\cal L}$ defines
in a unique way a minimal section $C = C(s)$
of the ruled surface
$S = {\bf P}_{C}(E)$ (see e.g. [H, ch.5, 2.6, 2.8]).
If $h^{0}(E) = 1$, the bundle section $s \in H^{0}(E)$
is unique, and the corresponding minimal section
$C(s)$ is unique. In contrary, if $h^{0}(E) \ge 2$,
the map

${\bf P}(H^{0}(E)) \rightarrow \{ \mbox{the minimal sections
of } S $, $s \mapsto C(s) \}$,

defines a nontrivial linear system of minimal sections of
$S$ (e.g., if $S$ is a quadric). Therefore, the set of
minimal sections of $S$ is the same as the set of the
bundle sections of normalized bundles which represent $S$.
In fact, if $g(C) \ge 1$ and $S$ is general, then
$h^{0}(E) = 1$ for any normalized $E$ which represents $S$.
In this case the questions (*) and (*) are equivalent.
This is exactly the setting of Theorem (2.4) for $k \ge 3$.

{\sc Definition.}
Call the line subbundle ( = the invertible subsheaf)
${\cal M} \subset E$
a {\it maximal subbundle} of $E$, if ${\cal M}$ is
a line subbundle of $E$ of a maximal degree.

\smallskip

Let $E$ be a fixed normalized bundle which represents $S$,
and let ${\cal M} \subset E$ be a maximal subbundle
of $E$.  Clearly  $deg({\cal M}) \ge 0$, since
${\cal O}_{C} \subset E$. Assume that $deg({\cal M}) >0$.
Then, after tensoring
by  ${\cal M}^{-1}$, we obtain the embedding
${\cal O}_{C} \subset E \otimes {\cal M}^{-1}$.

In particular, $h^{0}(E \otimes {\cal M}^{-1}) \ge 0$,
$E \otimes {\cal M}^{-1}$  represents $S$,  and
$deg(E \otimes {\cal M}^{-1}) < deg(E)$.  However
$E$ is normalized, hence $deg(E \otimes {\cal M}^{-1})$
cannot be less than $deg(E)$ -- contradiction. Therefore
$deg({\cal M}) = 0$,  and  the maximal subbundle ${\cal M}$
of $E$  defines the normalized bundle
$E \otimes {\cal M}^{-1}$  which also represents  $S$.

Therefore, we can reduce the question (*) to the
following question:

{\bf (***).} {\it How many maximal subbundles has a fixed
normalized rank 2 bundle $E$ which represents a given
ruled surface $S$ ?}

The answer of (*) - (***) for $S$ -- decomposable,
is given in [H, ch.5, Examples 2.11.1, 2.11.2, 2.11.3].
In particular, this implies the well known description
of the set of minimal sections of a rational ruled surface
$p : S \rightarrow {\bf P}^{1}$. In particular,
if $d = deg({\Delta}) = 4 \ or \ 5$, Theorem (2.4) implies that
the ruled surface $S(L)$ is either a quadric
(if $L \in C_{+}$), or the surface ${\bf F}_{1}$
(if $L \in C_{-}$). If $d = 6$,  the minimal sections
$C \in {\cal C}_{+} \cup {\cal C}_{-}$ are elliptic
curves.

If $d = deg({\Delta}) = 6$, the general
$C \in {\cal C}_{+} \cup {\cal C}_{-}$ is a smooth
elliptic curve. It follows from the proof of
Theorem (2.4) that $S(L)$ is a general element
of a versal deformation of ruled surface over an elliptic
base.  Similar conclusion takes place
if $deg({\Delta}) \ge 7$.  Therefore, the description
of the general fibers of ${\Phi}_{+}$ and ${\Phi}_{-}$
will follow from the description of the set of minimal
sections of the general member in a versal deformation
of a ruled surface.

\bigskip

{\bf (5.3)}
{\bf Lemma} (see [Se, Theorem 5]).
{\sl
Let $S \rightarrow C$ and $S' \rightarrow C'$ be two ruled
surfaces. Then $S$ and $S'$ can be deformed into each other
{\it iff} $C$ and $C'$ have the same genus, and the invariants
$e(S)$ and $e(S')$ have the same parity.
}

In particular, the parity of $e(S)$ is an invariant
of the deformations of the ruled surface $S$.

\bigskip

{\bf (5.4)}
{\bf Lemma} ( see [Se, Theorem 13]).
{\sl
The general surface in the versal deformation of a
rational ruled surface is a quadric if $e$ is
even, and the surface ${\bf F}_{1}$ if $e$ is odd.

The general surface of a versal deformation of a
ruled surface over elliptic base is a surface
represented by the unique indecomposable rank 2 vector
bundle of degree 1 if $e$ is odd, and a decomposable
ruled surface represented by a sum of two
(non-incident) line bundles of degree $0$ if $e$ is even.

The general surface of a versal deformation of a
ruled surface over a curve of genus $g \ge 2$ is
indecomposable. The invariant of such $S$ is
$g-1$ if $e \equiv g \mbox{ mod 2}$, or
$g$ if $e \equiv g-1 \mbox{ mod 2}$.
}

\bigskip

{\bf (5.5)}
{\bf Lemma} (see [H, ch.5, Example 2.11.2 and Exer.2.7]).
{\sl
Let $C$ be an elliptic curve, and let
let $S$ be the unique indecomposable ruled surface
over $C$ with invariant $e(S) = 0$. Then
the set ${\cal C}_{+}(S)$ of minimal sections
of $S$ form a 1-dimensional family parameterized
by the points of the base $C$. In particular, all
the minimal sections of $S$ are linearly
non equivalent.

Let $C$ be an elliptic curve, and let the ruled
surface $S$ be represented by the normalized
bundle $E = {\cal O}_{C} \oplus {\cal L}$,
where $deg({\cal L}) = 0$  and
${\cal L} \ne {\cal O}_{C}$. Then $S$ has
exactly two minimal sections: the section
$C = C(s_{E})$ defined by the unique bundle section
$s_{E}$ of $E$, and the section
$\overline{C}$ defined by the unique section
$s_{\overline{E}}$ of the second normalized bundle
$\overline{E} = {\cal O}_{C} \oplus {\cal L}^{-1}$
which represents $S$.
}

\bigskip

{\bf (5.6)}
{\bf Lemma} (see [LN, Proposition 2.4]).
{\sl
Let $S$ be an indecomposable ruled surface over
a curve $C$ of genus $g \ge 2$. Let $E$ be a fixed
normalized rank 2 bundle over $C$ which represents $S$,
and let
$[E] \in {\bf P}(H^{0}(K_{C} \otimes {\cal L}))$
be the point which corresponds to the extension
$0 \rightarrow {\cal O}_{C} \rightarrow E
\rightarrow {\cal L} \rightarrow 0$
defined by
$E$.
Let
${\alpha} :C \rightarrow {\bf P}(H^{0}(K_{C} \otimes {\cal L}))$
be the map defined by the linear system
$\mid K_{C} \otimes {\cal L} \mid$,
and let ${\alpha}(C)$ be the image of $C$.
Then
the set of maximal line subbundles ${\cal M}$
of $E$, which are different from ${\cal O}_{C}$,
is naturally isomorphic to the set of $e$-secant
line bundles of ${\alpha}(C)$ which pass through
the point $[E]$.
}

\bigskip

{\bf Remarks.}

{\bf  (1).}
Here $e = e(S) = deg({\cal L})$ is the invariant of $S$.

{\bf (2).}
The line bundle ${\cal D}$ on $C$ of degree $e$ is
an $e$-secant line bundle of ${\alpha}(C)$ which
passes through $[E]$, if the linear system
$\mid {\cal D} \mid$ contains an effective divisor $D$
such that the space $Span({\alpha}(D))$ passes through
the point [E].

\bigskip

{\bf (5.7)}
{\it The general position of the ruled surfaces $S_{L}.$}

\smallskip

Let $p : X \rightarrow {\bf P}^{2}$ be a smooth conic
bundle which induces a nontrivial unbranched double
covering ${\pi} : \tilde{\Delta} \rightarrow {\Delta}$
of the smooth discriminant curve $\Delta$ of degree
$deg({\Delta}) = d$. Without any restriction, we can
assume that $d \ge 4$ (otherwise the jacobian $J(X)$
will be trivial).
Let $C_{0} \subset {\bf P}^{2}$ be a general smooth
plane curve of degree $d-3$, and let
$S_{C_{0}} = p^{-1}(C_{0}) \subset X$ be, as usual,
the induced conic bundle surface over $C_{0}$.
In particular, we can assume that $C_{0}$ is smooth,
and $C_{0}$ intersects ${\Delta}$ in
$2g-2 = d(d-3)$ disjoint points
$(x_{1},...,x_{2g-2})$. Let
$L_{0} = x_{1}+...+x_{2g-2}$, and let
$F(L_{0}) = (S^{2g-2}{\pi})^{-1}(L_{0})
\subset S^{2g-2}(\tilde{\Delta})$.

As it follows from Theorem (2.4), the
assumption (2.4)(*) and the general choice of $C_{0}$
imply that the set $F(L_{0})$ splits into two disjoint
sets $F(L_{0})_{+}$ and $F(L_{0})_{-}$, of equal
cardinality = $2^{2g-3}$, with the following properties:

  {\bf 1.} If $L = l_{1} +..+ l_{2g-2} \in F(L_{0})_{-}$
then the ruled surface $S_{L}$ (obtained from $S_{C_{0}}$
by blowing down the $2g-2$ complimentary fibers
$\overline{l_{1}} ,..., \overline{l_{2g-2}}$) has invariant
$e = e(L) := e(S_{L}) = g-1$.

  {\bf 2.} If $L = l_{1} +..+ l_{2g-2} \in F(L_{0})_{+}$
then the ruled surface $S_{L}$ (obtained from $S_{C_{0}}$
by blowing down the $2g-2$ complimentary fibers
$\overline{l_{1}} ,..., \overline{l_{2g-2}}$) has invariant
$e = e(L) := e(S_{L}) = g$.

Moreover, as it follows from the
proof of (2.4), any such $S_{L}$ has to be a general element
of a versal deformation of a ruled surface over plane curve
of degree $d-3$, provided $C_{0}$ is chosen to be
sufficiently general. In particular, we can assume:

\bigskip

{\bf (5.7.1).}
(If $d \ge 7$, then)

{\sl

The general $S(L)$ is indecomposable; and if
$L \in F(C_{0})_{-}$ (see above), then:

{\bf (1.i).}
$e(S(L)) = g-1$;

{\bf (1.ii).} Let $E$ be a general normalized rank 2 bundle
which represents $S_{L}$.  Let
$s : 0 \rightarrow {\cal O}_{C_{0}} \rightarrow E$ be the
embedding defined by the bundle section $s$ of $E$,
and let ${\cal L}(s)$ be the cokernel of $s$.
Then the point
$[E] \in {\bf P}(H^{0}(K_{C_{0}} \otimes {\cal L}(s)))$
(which represents the extension class of $s$) is
in general position with respect to the set
of $(g-1)$-secant line bundles of ${\alpha}(C_{0})$
(see (5.6)).
}

\bigskip

{\it Remark.}

The condition (ii) is open, on the set of extension
classes $[E]$. Therefore, it has to be fulfilled for at least
one $C_{0} \in \mid {\cal O}_{{\bf P}^{2}} (d-3) \mid$.

Let e.g. $C_{0} \in {\bf S}_{d-3}$ (see (3.3)). The proof of
Lemma (3.3.1) implies that if such an $C_{0}$ is general, then
the centers ${\xi}_{ik}$ of the elementary transformations
$elm({\xi}_{ik})$  which interchange the ruled surfaces
$S_{L} : L \in F(C_{0})$ can be chosen to be general
(provided (2,4)(*) takes place).  We can assume that
$C_{0} \in {\bf S}_{d-3}$ has reduced components, which
intersect $\Delta$ transversely.  Now, it is not hard
to complete the definition of (normalized) bundles $E$,
bundle sections $s$,  extension classes $[E]$, etc.,
for such a curve $C_{0}$. In particular, the possibility
of the general choice of the centers ${\xi}_{ik}$
implies that (ii) takes place for such a curve $C_{0}$.
Therefore, the same is true also for the general
$C_{0} \in \mid {{\cal O}_{{\bf P}^{2}}} (d-3) \mid$.

In a more exact setting, if $(d-3) \le 3$,
the genus of $C_{0}$ is $0$ or $1$; and (i) and (ii)
have to be reformulated (in the obvious way)
in the context rational ruled surfaces  and ruled
surfaces over elliptic base (see e.g. (5.5)).

\smallskip

Similarly,

{\bf (5.7.2)}
Let $L \in F(C_{0})_{+}$, and let $C_{0}$ be general.
Then
the following takes place:

(If $d \ge 7$, then)

{\sl
{\bf (2.i).} $S_{L}$ is indecomposable, and
$e(S(L)) = g$;

{\bf (2.ii).}
The point which is defined by the general
extension class $[E]$ is in general
position with respect to the set of
$g$-secant line bundles of ${\alpha}(C_{0})$
}
(see (5.7.1)(1.ii)).

\bigskip

{\bf (5.7.3).}
For $d \le 6$ -- see the last remark.

\smallskip

Now, (5.1) - (5.7) and Theorems (2.4) and (4.4)
imply the following

\bigskip

{\bf (5.8).}
{\it Description of the general fibers of the
Abel-Jacobi maps ${\Phi}_{+}$ and ${\Phi}_{-}$.}

\smallskip

{\bf Theorem.}
{\sl
Let $p : X \rightarrow {\bf P}^{2}$ be a conic bundle,
as in (5.1).
Let
${\cal C}_{+}$ and ${\cal C}_{-}$ be the families of
non-isolated and isolated minimal sections, and let
${\phi} : {\cal C}_{+} \rightarrow C_{+}$,
${\phi} : {\cal C}_{-} \rightarrow C_{-}$,
${\psi} : Supp({\Theta}) \rightarrow {\Theta}$, and
${\psi} : Supp(P^{-}) \rightarrow P^{-}$
be the families and the natural maps defined in (4.1).
Let ${\Phi}_{+} : {\cal C}_{+} \rightarrow J(X)$
and ${\Phi}_{-} : {\cal C}_{-} \rightarrow J(X)$
be the Abel-Jacobi maps for ${\cal C}_{+}$
and ${\cal C}_{-}$, and let $Z_{+}$ and $Z_{-}$
be the images of ${\Phi}_{+}$ and ${\Phi}_{-}$.

Then one of the following two alternatives is true:

  {\bf (A,+).}
$C_{+} = Supp({\Theta})$,
$Z_{+}$ is a translate of $\Theta$
( $\Leftrightarrow C_{-} = Supp(P^{-})$,
$Z_{-} = J(X) \cong P$).

Let
$z \in Z_{+}$ is general, and let
${\cal L} = j \circ i^{-1} (z) \in {\Theta}$ be the
sheaf which corresponds to $z$ ( see the proof of (4.4)).
Then:

  (1). The fiber
${\cal C}_{+}(z) := {\Phi}_{+}^{-1}(z)$ is 2-dimensional.

  (2.) The map $\psi$ defines on
${\cal C}_{+}(z)$ the natural fibration
${\psi}(z) : {\cal C}_{+}(z) \rightarrow \mid {\cal L} \mid
\cong {\bf P}^{1}$.

  (3). The general fiber
${\cal C}_{+}(L) := {\psi}(z)^{-1}(L)$ of ${\psi}(z)$
can be described as follows ($d \ge 4$):

Let $C_{0}(L) \subset {\bf P}^{2}$ be the plane
curve of degree $d-3$ defined by $L$. Then

  (i). If $d = deg({\Delta}) = 4 \mbox{ or } 5$, then
$S(L) \cong {\bf P}^{1} \times {\bf P}^{1}$,  and
${\cal C}_{+}(L) \cong$ the fiber ${\bf P}^{1}$ of
the projection
$p(L) : S(L) \rightarrow C_{0}(L) \cong {\bf P}^{1}$
induced by $p$;

  (ii). If $d = deg({\Delta}) = 6$, then
$p(L) : S(L) \rightarrow C_{0}(L)$ is the only
indecomposable ruled surface over the elliptic
case $C_{0}(L)$, and the fiber
${\cal C}_{+}(L)$
of
${\psi}(z) : {\cal C}_{+}(z) \rightarrow
\mid {\cal L} \mid \cong {\bf P}^{1}$
is isomorphic to $C_{0}(L)$.
In particular, ${\cal C}_{+}(z)$
is an elliptic fibration over the rational base
curve $\mid {\cal L} \mid$;

  (iii). Let $d = deg({\Delta}) \ge 7$,  let
$g = d(d-3)/2 + 1$ be the genus of $C_{0}(L)$,
let $C \in {\cal C}_{+}(L)$ be general, and let
$0 \rightarrow {\cal O}_{C_{0}} \rightarrow E
\rightarrow {\cal N} \rightarrow  0$
be the extension defined by the section $C$.  Let
${\alpha}(C_{0}) \subset
{\bf P}(H^{0}(K_{C_{0}} \otimes {\cal N}))$
be the image of $C_{0}$
defined by the sheaf $K_{C_{0}} \otimes {\cal N}$.

Then ${\bf P}(H^{0}(K_{C_{0}} \otimes {\cal N}))$
$\cong {\bf P}^{2g-2}$,  ${\alpha}$
is a regular morphism of degree 1, and the
point $[E]$ defined by this extension is in
general position with respect to the set
of $g$-secant line bundles of ${\alpha}(C_{0})$.
Moreover,
${\cal C}_{+}(L)$ is birational to the
1-dimensional set $\{ g-sec({\alpha},E) \}$  of $g$-secant planes
of ${\alpha}(C_{0})$ through the
point $[E]$. In particular, if $C'$ and $[E']$
is another pair of this type, then the
normalizations of the curves  $\{ g-sec({\alpha},E) \}$ and
$\{ g-sec({\alpha}', E') \}$ are isomorphic to each other.

  (4). If $z \in Z_{-}$ is general and
${\cal L} = j \circ i^{-1}(z)$, then
$\mid {\cal L} \mid \cong{\bf P}^{0}$.
If $L = L(z)$ is the unique element
of $\mid {\cal L} \mid$, then
the fiber
${\Phi}_{-}(z) = {\Phi}_{-}^{-1}(z)$
is discrete and:

  (i). If $d = deg({\Delta}) = 4 \mbox{ or } 5$, then
${\Phi}_{-}(z)$ has exactly one element --
defined by the unique $(-1)$-section of
the ruled surface $S(L) \cong {\bf F}_{1}$.

  (ii). If $d = deg({\Delta}) = 6$, then
${\Phi}_{-}(z)$ has exactly two elements --
defined by the two (disjoint) sections
of the decomposable ruled surface
$S_{L}$ over the elliptic base $C_{0}(L)$.

  (iii). Let $d = deg({\Delta}) \ge 7$.
Then the fiber ${\Phi}_{-}(z)$ is
isomorphic to the fiber
${\psi}^{-1}(L)$.
Let $C$ be some element of this
fiber, let
$0  \rightarrow {\cal O}_{C_{0}(L)} \rightarrow E
\rightarrow {\cal N} \rightarrow 0$
be the extension defined by $C$, and let
${\alpha}(C_{0})$ and $[E]$ be as above.
Then $\mid K_{C} \otimes {\cal N} \mid$
$\cong {\bf P}^{2g-3}$, ${\alpha}$ is
generically of degree 1, and $[E]$ does not
lie on an infinite set of $(g-1)$-secant
planes of ${\alpha}(C_{0})$.
Moreover,  the cardinality
of  ${\cal C}_{-}(z)$ is equal to
$\# \{ (g-1)- \mbox{secant planes of } {\alpha}(C_{0}) \} + 1$
(see (5.6)).

{\bf (A,--)}.
$C_{-} = Supp({\Theta})$,
$Z_{-} \cong {\Theta}$
($\Leftrightarrow$
$C_{+} = Supp(P^{-})$,
$Z_{-} = J(X) \cong P$).

Then the description of the general fibers of
${\Phi}_{-}$ and ${\Phi}_{+}$ is similar to
this from (A.+)(1)-(4). We shall mark only the
differences:

(1)-(2)-(3). The fiber ${\cal C}_{-}(z)$ is 1-dimensional.
The map
${\psi}(z) : {\cal C}_{-}(z) \rightarrow \mid {\cal L}(z) \mid$
$\cong {\bf P}^{1}$ is finite and surjective,
and the fiber of ${\psi}(z)$ has the same description
as the fiber
${\cal C}_{-}(L) = {\psi}^{-1}(L)$ described
in (A.+)(4).

(4). The fiber ${\cal C}_{+}(z)$ is 1-dimensional.
Let ${\cal L} = {\cal L}(z) = j \circ i^{-1}(z)$, and
let $L = L(z)$ be the unique element of the linear
system $\mid {\cal L} \mid$.  Then the sets
${\cal C}_{+}(z)$ and ${\cal C}_{+}(L)$ coincide.
In particular, the fiber ${\cal C}_{+}(z)$ has
the same description as the fiber
${\cal C}_{+}(L)$ described in (A.+)(1)-(4).
Note that here $L = L(z)$ is unique.
}

\bigskip

{\bf (5.9)}
{\bf Corollary.}
{\sl
$dim ({\cal C}_{+}) = dim \ J(X) + 1$ =
$dim \ ( \mid {{\cal O}_{{\bf P}^{2}}} (d-3) \mid ) + 1 =
(d-2)(d-1)/2 .$

$dim({\cal C}_{-}) = dim \ J(X)$ =
$dim \ ( \mid {{\cal O}_{{\bf P}^{2}}} (d-3) \mid ) =
(d-2)(d-1)/2 - 1 .$
}

%%%%%%%%%%%%%%%%%%%%%%%%%%%%%%%%%%%%

\newpage

%%%%%%%%%%%%%%%%%%%%%%%%%%%%%%%%%%%%

\centerline{\bf 6. Examples.}

\bigskip
\bigskip

{\bf (6.1)}
{\sc
The double covering of ${\bf P}^{2} \times {\bf P}^{1}$.
}

\smallskip

{\bf (6.1.1)}
{\sc Definition.}
Let $Y = {\bf P}^{2} \times {\bf P}^{1} \subset {\bf P}^{5}$,
let $Q$ be a general smooth quadric in ${\bf P}^{5}$, and let
${\xi} : X \rightarrow Y$ be the double covering branched along
the smooth surface $S = Y \cap Q$.

Clearly, $X$ is a Fano 3-fold with $rank \ {\bf Pic}(X) =2.$
We call the  threefold $X$  {\it the double covering of}
${\bf P}^{2} \times {\bf P}^{1}.$

\smallskip

{\bf (6.1.2)}
{\it The conic bundle structure on $X,$
and the families ${\cal C}^{0}_{1,0}$ and ${\cal C}^{0}_{1,1}.$}

\smallskip

Let $p_{o} : Y \rightarrow {\bf P}^{2}$ and
$q_{o} : Y \rightarrow {\bf P}^{1}$ be the projections, and
let $p = p_{o} \circ {\xi}$, $q = q_{o} \circ {\xi}$.

Then $p : X \rightarrow {\bf P}^{2}$ is a conic bundle,
and the discriminant ${\Delta}$ is a smooth plane curve
of degree $4$. Let ${\eta} \in {\bf Pic}^{0}_{2}({\Delta})$
be the torsion sheaf which defines the double covering
${\pi} : \tilde{\Delta} \rightarrow {\Delta}$. Then the
Prym-canonical system $\mid K_{\Delta} + {\eta} \mid$ is
naturally isomorphic to ${\bf P}^{1}$. Indeed, the elements
of the system
$\mid K_{\Delta} + {\eta} \mid$ =
$\mid {\cal O}_{\Delta} (1 + {\eta}) \mid$
are in (1:1)-correspondence with the ${\bf P}^{1}$-family
of conics $q$ which are totally tangent to $\Delta$ along
effective divisors $\delta$ of the linear system
$\mid {\cal O}_{\Delta} (1 + {\eta}) \mid$. The family of all
these conics is parameterized naturally by the points
$x \in {\bf P}^{1}$. Indeed, let $x \in {\bf P}^{1}$ be general,
and let $Q_{x} = q^{-1}(x)$. Then $Q_{x}$ is a quadric, and
the projection $p$ defines a double covering
$q : Q_{x} \rightarrow {\bf P}^{2}$ branched along a conic
$q(x)$ which is totally tangent to $\Delta$ along
a degree $4$ divisor ${\delta}(x)$.  It is
not hard to see that
${\delta}(x) \in \mid {\cal O}_{\Delta}(1 + {\eta}) \mid$.
In particular, the pencil of quadrics
$q : X \rightarrow {\bf P}^{1}$ has 6 degenerations -
the cones $Q_{x_{i}} , i = 1,...,6$ which correspond to
the $6$ degenerated conics $q(x_{i})$ in the quadratic
pencil $\{ q(x) \}$ of conics totally tangent to the plane
quartic $\Delta$.
Let $x \in {\bf P}^{1}$ be general, and let
$p : Q_{x} \rightarrow {\bf P}^{2}$ be the double covering
branched along the smooth conic $q(x)$.
The double covering $p$ maps the elements $C$ of the two
generators ${\Lambda} (x)$ and $\overline{\Lambda} (x)$
of $Q_{x}$ to lines $l = p(C)$ tangent to $q(x)$. Clearly,
these elements $C$ are (1,0)-lines on $X$, i.e. - curves
of bidegree (1,0) on $X$ (with respect to the bidegree
map induced by the projections $p$ and $q$). Moreover, any
such $C$ is a section of the conic bundle surface
$S_{l} = p^{-1}(l)$, where $l= l(C)$. It follows that
the 2-dimensional family ${\cal C}^{0}_{1,0}(X)$ of
(1,0)-lines on $X$ coincides with the family
${\cal F}_{-}  = {\cal C}_{-}$ of isolated sections
over the canonical system of $\Delta$. The including
of any (1,0)-line $C$ in a generator of a quadric $Q_{x}$
implies that the effective divisor
$L(C)$ on $\tilde{\Delta}$ defined by
${\psi}(C) = C \cap p^{-1}({\Delta})$ moves in a non-trivial
linear system of dimension one. Therefore
${\psi}({\cal C}^{0}_{0,1}(X)) = Supp({\Theta})$, and
the set of fibers of ${\psi}$ is in (1:1)-correspondence with
the set of generators $ Gen (X)$ = (the closure of)
$\{ ({\Lambda}(x), \tilde{\Lambda}(x)):$
$x \in {\bf P}^{1} - \{ x_{1},...,x_{6} \} \}$.
Clearly, these generators can be treated also as
sheaves on $\tilde{\Delta}$ which belong to ${\Theta}$.
Moreover, since $dim(J(X)) = 2$,
$J(X)$ is an jacobian of a curve $C_{X}$ of genus $2$.
In particular $C_{X} \cong {\Theta} \cong  Gen(X)$,
and the natural double covering
$q : Gen(X) \rightarrow {\bf P}^{1}$ defines
a double covering $q' : C_{X} \rightarrow {\bf P}^{1}$,
branched in the $6$ points $x_{1},..,x_{6}$ over which
the $({\bf P}^{1} \times {\bf P}^{1})$-bundle
$q : X \rightarrow {\bf P}^{1}$ degenerate.

Now, it is not hard to see that:

\smallskip
{\bf (6.1.3)}
{\sl
The 3-dimensional family
${\cal C}_{+}$ of non-isolated sections (over the canonical
system of $\Delta$) coincides with the family
${\cal C}^{0}_{1,1} (X)$ of bidegree $(1,1)$-conics on $X$.
In particular, the Abel-Jacobi map
${\Phi}_{+} : {\cal C}^{0}_{1,0} (X) \rightarrow J(X)$
is surjective, and the general fiber of ${\Phi}_{+}$
is isomorphic to ${\bf P}^{1}$.
}

\bigskip
\bigskip

{\bf (6.2)}
{\sc
The intersection of two quadrics.
}

\smallskip
{\bf (6.2.1)}
Let $X = X_{2.2} \subset {\bf P}^{5}$ be a smooth intersection
of two quadrics.
The threefold $X$ is a Fano 3-fold with
$rank \ {\bf Pic}(X) = 1.$
In particular, there are no biregular conic bundle structures on
$X$ (otherwise, $rank \ {\bf Pic}(X) \ge 2).$
However, $X$ admits various natural birational conic bundle
structures indexed by the general elements of the 4-dimensional
family of conics on $X.$

\smallskip
{\bf (6.2.2)}
{\it The conic bundle structure defined by a conic $q_{o}$
on $X$, and the families ${\cal C}^{0}_{1}(X)$ and
${\cal C}^{0}_{2}(X)_{{q}_{o}}$}.

\smallskip

Let $q_{o}$ be a sufficiently general fixed conic on
$X = X_{2.2} \subset {\bf P}^{5}$, and let
$p_{o} : X \rightarrow {\bf P}^{2}$ be the rational
projection centered in the plane of the conic $q_{o}$.
Then $p_{o}$ defines a birational structure of
a conic bundle $p : X \rightarrow {\bf P}^{2}$,
and the non-trivial component $\Delta$ of the determinant
locus of $p$ is a smooth plane quartic.

The canonical system of $\Delta$ is defined by the system
of lines $l$ in ${\bf P}^{2}$; and it is not hard to see that
the ``canonical'' family
${\cal S}[1]$ of conic bundle surfaces over the lines in
${\bf P}^{2}$ is versally embedded.
Let ${\cal C}_{-} = {\cal F}_{-}$
and ${\cal C}_{+} = {\cal F}_{+}$ be the families
of isolated and non-isolated minimal sections of the elements
of the system ${\cal S}[1]$.  By construction,  the general
curve $C \ \in {\cal C}_{-} \cup {\cal C}_{+}$ is projected
isomorphically (by $p_{o}$) onto a line
$l \subset {\bf P}^{2}$. The same property have the
general elements $C$ of the families

${\cal C}^{0}_{1} (X)$ =
$\{ \mbox{ the lines in } X \}$  and
${\cal C}^{0}_{2} (X)_{q_{o}}$ =
$\{ \mbox{ the conics in } X \mbox{ which intersect } q_{o} \} .$

Moreover,  $dim({\cal C}_{+}) = 3 =
dim({\cal C}^{0}_{2} (X)_{q_{o}})$,  and
$dim({\cal C}_{-}) = 2 =
dim({\cal C}^{0}_{1} (X))$.

Therefore, by disregarding the necessary corrections
caused by the non-standard form of the conic bundle
projection $p$, we can assume that
${\cal C}_{+} \cong {\cal C}^{0}_{2}(X)_{q_{o}}$,  and
${\cal C}_{-} \cong {\cal C}^{0}_{1}(X)$.

\smallskip

{\bf (6.2.3)}
{\it The parameterization of $(J(X),{\Theta})$
via minimal sections.}

\smallskip

It is well known ( see e.g. [Do]) that
the family ${\cal C}^{0}_{1}(X)$ is isomorphic
to the intermediate jacobian $J(X)$. Therefore,
${\Phi}_{-}$ sends the 3-dimensional family
${\cal C}_{+}$ onto a copy of the theta divisor
${\Theta}(X)$, i.e. (in this case) the theta divisor
is described as an Abel-Jacobi image of the
family of non-isolated minimal sections of the canonical
system of conic bundle surfaces ${\cal S}[1]$.

In fact, the conics on $X = X_{2.2}$ are described
by the elements of the generators of the
${\bf P}^{1}$-system of quadrics in ${\bf P}^{5}$
which contain $X$. The general of these quadrics
is a quadric $Q(x)$ of rank $6$, and $Q(x)$ has two
generators
${\Lambda}(x) \cong \overline{\Lambda}(x)
\cong {\bf P}^{2}$. Any of these two generators
define a ${\bf P}^{2}$-system of conics
in $X$.  In particular, any conic on $X$ is rationally
equivalent to a conic on $X$ which intersects the
fixed conic $q_{o}$.

Therefore the Abel-Jacobi
image of the 4-dimensional family
${\cal C}^{0}_{2}(X)$ of conics on $X$ isomorphic to
the Abel-Jacobi image of ${\cal C}^{0}_{2} (X)_{q_{o}}$.
Moreover, the induced (by the generators) double
covering of ${\bf P}^{1}$ is branched over the set
of $6$ points which represent the $6$ singular quadrics
through $X$. Just as in (6.1), the so defined curve
${\ Gen}$ is isomorphic to $\Theta$, and (since
$dim(J(X)) = 2$) $J(X) \cong J({\ Gen})$ as p.p.a.v.

\bigskip
\bigskip

{\bf (6.3)}
{\sc
The cubic threefold.
}

\smallskip

{\bf (6.3.1)}
Let $X = X_{3} \subset {\bf P}^{4}$ be a general smooth
cubic threefold, and let ${\cal F}$ be the 2-dimensional
family of lines on $X$. It is well known
(see e.g. [CG] or [Be2])
that the Abel-Jacobi map

${\Phi}^{+,-} : {\cal F} \times {\cal F} \rightarrow J(X)$,
${\Phi}^{+,-}: (l , m) \mapsto {\Phi}([l - m])$

is a 6-sheeted covering of a copy of the theta divisor
${\Theta}$ of $J(X)$, and the general fiber of
${\Phi}^{+,-}$ is a Schl\"afli's double-six of lines
on a hyperplane section of $X$.

We shall see (by using appropriate conic bundle
structures on $X$) that the Abel-Jacobi image of the
6-dimensional family ${\cal C}^{0}_{3}(X)$
of twisted cubics on $X$ is also a copy of the
theta divisor $\Theta$ of the intermediate jacobian $J(X)$.
The advantage of the second parameterization of $\Theta$
is that the fibers of the Abel-Jacobi map
${\Phi} : {\cal C}^{0}_{3} \rightarrow {\Theta}$ are
connected (the general one is, in fact, a projective plane --
see below).

\newpage

{\bf (6.3.2)}
{\it The families ${\cal C}_{+}$ and ${\cal C}_{-}$ for the
birational conic bundle structure defined by a line
$l_{0} \subset X.$}

\smallskip

Let $l_{0}$ be a sufficiently general fixed line on $X$,
and let $p_{o} : X \rightarrow {\bf P}^{2}$ be the
rational projection from the line $l_{0}$.  Let
$\tilde{X}$ be the blow-up of $X$ along $l_{0}$.
The map $p_{o}$ defines a conic bundle structure
$p : \tilde{X} \rightarrow {\bf P}^{2}$
on $\tilde{X}$, and the nontrivial component
$\Delta$ of the discriminant curve of $p$
is a smooth plane quintic.

In this case, the families ${\cal C}_{-}$ and
${\cal C}_{+}$ are the families of isolated
and non-isolated sections of the 5-dimensional family
${\cal S}[2]$ of conic bundle surfaces over the
${\bf P}^{5}$-space of plane conics.
Just as in (6.2), we see that (outside the components
and subsets caused by the blow up and some coincidences
with $l_{0}$):

\smallskip

  (1). The family
${\cal C}_{+}$ can be identified with the 6-dimensional
family ${\cal C}^{0}_{4} (X)_{l_{0}, l_{0}}$ of
rational quartics on $X$ which intersect twice the line
$l_{0}$.

\smallskip

  (2).  The family
${\cal C}_{-}$ can be identified with the 5-dimensional
family ${\cal C}^{0}_{3} (X)_{l_{0}}$ of
twisted cubics on $X$ which intersect the line $l_{0}$.

\smallskip

Let $C \subset X$ be a general twisted cubic on $X$
which intersects the line $l_{0}$,  and let
$x_{0}$ be the common point of $C$ and $l_{0}$.
Let $H = Span (C)$ be the hyperplane spanned on $C$,
and let $S_{H} = X \cap H$ be the hyperplane section
defined by $H$. The general choice of $C$ implies that
$S_{H}$ is a smooth cubic hypersurface. Since $C$
is a twisted cubic on $S_{H}$, the curve $C$ moves
in a $2$-dimensional linear system on $S_{H}$.
Again by the general choice of $C$, the non-complete
linear system
$\mid { {\cal O}_{S_{H}} } (C) - x \mid$
is isomorphic to ${\bf P}^{1}$.

Therefore  the natural ``intersection'' map :

${\psi} : {\cal C}_{-} \cup {\cal C}_{+}
\rightarrow Supp({\Theta}) \cup Supp(P^{-})$
sends the general element $C \in {\cal C}_{-}$
onto an element $L(C) = {\psi}(C)$ of a
non-trivial linear system.

Therefore  ${\Phi}_{+}({\cal C}_{-}) = {\Theta}$,
${\Phi}_{+}({\cal C}_{+}) = J(X)$.

The same argument as in (6.2) (based on the existence
of a rational deformation to an element of a subfamily)
implies

\smallskip

{\bf (6.3.3). Proposition.}
{\sl
Let $X$ be a general smooth cubic hypersurface in
${\bf P}^{4}$.
Let ${\cal C}^{0}_{3} (X)$ be the 6-dimensional family
of twisted cubics on $X$, and let
${\cal C}^{0}_{4} (X)$ be the family of rational
quartics on $X$. Then

  (1). The Abel-Jacobi map ${\Phi}^{0}_{3}$ for
${\cal C}^{0}_{3}$
sends ${\cal C}^{0}_{3} (X)$ onto a copy of the theta
divisor $\Theta$ of $J(X)$.

Moreover, if $C \in {\cal C}^{0}_{3}$ is general, and
$z = {\Phi}^{0}_{3}(C)$, then the fiber
$({\Phi}^{0}_{3})^{-1}(z)$ coincides with
the  2-dimensional complete linear system
$\mid { {\cal O}_{S} } (C) \mid$ on the
smooth cubic surface $S = X \cap Span(C)$.

  (2). The Abel-Jacobi map ${\Phi}^{0}_{4}$ for
${\cal C}^{0}_{4}$
sends the family ${\cal C}^{0}_{4} (X)$ onto the
intermediate jacobian $J(X)$.
}

\smallskip

{\bf (6.3.4) Remarks.}

  {\bf (1).}
{\it The degree of the Gauss map for ${\Theta}$ via
the family of twisted cubics.}

\smallskip

Let $\gamma$ be rational Gauss map

${\gamma} : {\Theta} - Sing({\Theta}) \rightarrow
{\bf P}(T_{J(X)})_{0}^{*} \cong { {\bf P}^{4} } ^{*}$ ,
${\gamma} : z \mapsto
{\bf P}( \mbox{the tangent space of } {\Theta}
\mbox{ in } z )$,

and let ${\cal C}^{0}_{3} \rightarrow { {\bf P}^{4} }^{*}$
e the (rational) span-map defined, on the general
element $C$ of ${\cal C}^{0}_{3}$,
by:
$C \mapsto Span(C)$.
It can be seen (see e.g. [Vo], [T],
[CG])  that
$Span = {\gamma} \circ {\Phi}^{0}_{3}$.
In particular,
we obtain:

\smallskip

{\sl
The degree of the Gauss map
${\gamma}_{\Theta}$, for the theta divisor of
a cubic threefold, is $72$ = the number of
linear systems of twisted cubics on a smooth
cubic surface.
}

\smallskip

In connection with the parameterization of $\Theta$ by
$F \times F$ (see (6.3.1)),
it is not hard to find a (1:1)-correspondence between the
set of Schl\"afli's double-six'es on a fixed cubic
surface $S$, and the set of linear systems of twisted cubics
on $S$. In fact, the double-six'es on $S$ are in
(1:1)-correspondence with the set of morphisms
${\sigma} : S \rightarrow {\bf P}^{2}$. The linear system
$\{ C \} ({\sigma})$ of twisted cubics which corresponds
to such $\sigma$ is the $\sigma$-preimage of the system
$\mid { {\cal O}_{P^{2}} }(1) \mid$.

\newpage

  {\bf (2).}
{\it The point $Sing({\Theta})$}.

It is well-known (see e.g. [CG], [Tju]) that $\Theta$
has unique singular point $o$, and the base of the
tangent cone to $\Theta$ in $o$ is isomorphic to the
cubic $X.$  The point $o$
coincides with the image of the 5-dimensional subfamily
${\cal D} \subset {\cal C}^{0}_{3}$ of nodal plane cubics
on $X$. Indeed, ${\cal D}$ coincides with the set
of these $C \in {\cal C}^{0}_{3}$ where the rational Gauss
map ${\gamma}_{\Theta}$ is not regular on ${\Phi}(C)$
(see (1)), etc.

\bigskip
\bigskip

{\bf (6.4)}
{\sc
The Fano threefold $X_{16} \subset {\bf P}^{10}$.
}

\smallskip

{\bf (6.4.1)}
{\it The conic bundle structure on $X_{16}.$}

\smallskip

Let $W \subset {\bf P}^{10}$ be a general hyperplane section
of the Segre 5-fold
${\bf P}^{2} \times {\bf P}^{3} \subset {\bf P}^{11}$,
and let $X \subset W$ be a general smooth divisor
of bidegree $(1,2)$. Let $p : X \rightarrow {\bf P}^{2}$
and $q : X \rightarrow {\bf P}^{3}$ be the standard
projections.  It can be easily verified that the map $p$
defines
a standard conic bundle structure on $X$,
and the map $q$ is a blow up of a smooth curve
$B \subset {\bf P}^{3}$ of genus $5$, and of degree $7$.
Denote by $l$ and $h$ be the generators of
${\bf Pic}({\bf P}^{2})$ and ${\bf Pic}({\bf P}^{3})$
(and also their preimages on the various subsets of
${\bf P}^{2} \times {\bf P}^{3}$).

The invariants of $X$ can be computed in many ways, e.g.:

By adjunction
$- K_{X} = (l+h)$ is a hyperplane
section of $X$, i.e. $X$ is a Fano threefold with
$rank({\bf Pic} (X)) = 2$. The identities
$l^{3} = 0,$  $h^{4} = 0$,  $l^{2} h^{3} = 1$, and
$[X] = (\mbox{ the class of } X) = (l+h)(l+2h)$,
imply that $(-K_{X})^{3} = 16$.
The invariant $h^{2,1}(X) = 5$ can be computed from the
tangent bundle sequence for the embedding
$X \subset {\bf P}^{2} \times {\bf P}^{3}$.
In the list of Mori and Mukai (see [MM]), the only Fano
3-fold with these invariants is the blow up of ${\bf P}^{3}$
along a smooth curve $B$ of degree $7$ and of genus $6$,
s.t. $B \subset {\bf P}^{3}$ is a (non-complete)
intersection of cubics.

In fact, the blow up coincides with $q$, and
the map
$q \circ p^{-1}$, $l \mapsto q(S_{l}) = q(p^{-1}(l))
\subset {\bf P}^{3}$

defines an isomorphism between ${{\bf P}^{2}}^{*}$ and the
non-complete linear system
$\mid { {\cal O}_{P^{3}} } (3 - B) \mid$

of cubics through $B$. The degree $deg({\Delta})$ of the
discriminant curve ${\Delta}$ of the standard conic bundle
$p : X \rightarrow {\bf P}^{2}$ can be computed
by the
formula:
$-4.K_{{\bf P}^{2}} \equiv p_{*}[(-K_{X})^{2}] + {\Delta}$,

where $\equiv$ is a num.equivalence $=$ a lin.equivalence,
on ${\bf P}^{2}$ (see [MM]). Therefore

$deg({\Delta}) = 12 - l.(-K_{X})^{2}.[X] =
12 - l(l+h)^{2}(l+h)(l+2h) = 5$.
The discriminant curve $\Delta$ is smooth, since
$5 = h^{2,1}(X) = dim(J(X))$ =
$dim(P(\tilde{\Delta},{\Delta}))$ = $g({\Delta}) - 1$.

\smallskip

{\bf (6.4.2)}
{\it The families ${\cal C}_{+}$ and ${\cal C}_{-}.$}

\smallskip

Let $x \in {\bf P}^{2}$ be general, and let $f_{x}$ be the
conic $f_{x} = p^{-1}(x)$.  Let $l$ and $m$
be two lines such that $l \cap m = x$.
Then the curve $q(S_{l}).q(S_{m}) = B + q(f_{x})$ is a
complete intersection of two cubics.  The adjunction formulae
for the arithmetical genus of a complete intersection
imply that the conic $q(f_{x})$ intersects $B$ in $6$
points.  Now, it is not hard to see that the
${\bf P}^{2}$-system of $(0,2)$-conics on $X$ is mapped
isomorphically onto a system of $6$-secant conics
of $B$.  Similar computations imply that:

\bigskip

{\bf (1).}
{\sl
   $\tilde{\Delta}$ = ${\cal C}^{0}_{0,1}$ =
(the set of $(0,1)$-lines on $X$)
$\cong$ (the set of $3$-secant lines of $B$).

  The $5$-dimensional family ${\cal C}_{-}$ of isolated
minimal sections of the system ${\cal S}[2]$ coincides
with the family ${\cal C}^{0}_{2,3} (X)$ of
rational curves of bidegree $(2,3)$ on $X$,
and $q$ maps the last family on a component
${\cal C}^{0}_{3}[7](B)$ of twisted cubics which
are $7$-secant of $B$;

  Similarly, ${\cal C}_{+}$ = ${\cal C}^{0}_{2,4}$,
and $q$ defines an isomorphism between this family and a
component ${\cal C}^{0}_{4}[10](B)$ of
the family of rational quartic
curves which are $10$-secant of $B$.
}

\bigskip

In contrast with the example (6.3), the set
$Supp({\Theta})$ coincides with
$C_{+} = {\psi}({\cal C}_{+})$,  i.e.:

\smallskip

{\bf (2).}
{\sl
The theta divisor $\Theta$ of $J(X)$ is
parameterized by the Abel-Jacobi image of the
$6$-dimensional family of non-isolated
minimal sections of the system ${\cal S}[2]$.
}

\smallskip

{\it Proof of (2).}
The verification of the coincidence (2) is reduced
to the verification of the fact that
the general twisted cubic $C_{3}$,  which is $7$-secant
of $B$, cannot move in a rational system
in ${\cal C}^{0}_{3}[7](B)$.

Let, in contrary, $Supp({\Theta}) = C_{-}$, and let
$\{ C(t) : t \in {\bf P}^{1} \}$ be a general
(possibly -- non-linear) pencil
in ${\cal C}^{0}_{2,3}$ which represents a point
of $\Theta$. The $1$-dimensional family  $\{ C(t) \}$
has a finite number of degenerations of type
(section + section). It is easy to see that all these
degenerations have to be of the form
$C^{0}_{1,1} + C^{0}_{1,2}$,  where $C^{0}_{1,1}$
is a conic on $X$ of bidegree $(1,1)$ (which is mapped
onto a  bisecant line of $B$), and $C^{0}_{1,2}$
is a twisted cubic on $X$ of bidegree $(1,2)$ (which is
mapped onto a $5$-secant conic of $B$).
Any of these twisted cubics $C = C^{0}_{1,2}$ is
a component of the 1-cycle
$w(C) = p^{-1}(p(C)) \cap q^{-1}(q(C))$ on $X$.

On the one hand, the computation
in the Chow ring  $A.({\bf P}^{2} \times {\bf P}^{3})$
implies that the class
$[w(C)] \in A.({\bf P}^{2} \times {\bf P}^{3})$,
of the cycle $w(C)$,  is $4lh^{3} + 6l^{2}h^{2}$
(in $A.(X)$ this class is $2lh$).

On the other hand, $w(C) = C + \overline{C}$, where
$\overline{C}$ is a $(1,1)$-conic on $X$; in particular,
$\overline{C}$ is an isolated section of the system ${\cal S}[1]$,
i.e. $\overline{C} \in {\cal F}_{-}$ is an isolated Fano
section on $X$.
Therefore, $[C] = [\overline{C}] - [2lh]$.
By applying this to the general pencil in
${\cal C}_{-}$, and by exploring the fact that the
Abel-Jacobi maps are factorized through the cycle-class map,
we
obtain
that
${\Theta}$ =
${\Phi}_{+}({\cal C}^{0}_{2,3})$ is a translate of
${\Phi}^{+,-}({\cal F}_{-} \times {\cal F}_{-})$,

\smallskip

where
${\Phi}^{+,-} : (C_{1},C_{2}) \mapsto {\Phi}([C_{1} - C_{2}])$ ,
and $\Phi$ is the quotient Abel-Jacobi map for the
corresponding cycle-classes.
By using the fact that the Fano family
${\cal F}_{-}$ is mapped onto the family of bisecants of $B$,
and by the obvious isomorphism
$(J(X), {\Theta}) \cong (J(B), {\Theta}(B))$,
we obtain that ${\Theta}(B)$ is parameterized by
(the linear equivalence classes of the divisors of) the
$4$-dimensional
family
$\{ x_{1} + x_{2} - x_{3} - x_{4} : x_{i} \in B \}.$

However, this is impossible. Indeed, the initial general
nontrivial linear system
$\{ C(t) \} \subset {\cal C}_{-}$ is projected via $p$
to a non trivial system of plane conics
$\{ p(C(t)) \}$. This system forms a rational curve
$A$ of degree $d \ge 1$,  in the ${\bf P}^{5}$-space of
plane conics.  Therefore, the curve $A$ intersects the
determinantal cubic $Det$ of all degenerate plane
conics in $3d$ points. Therefore, if
$x_{i} , i = 1,...,4$ is a general $4$-tuple of points
of $B$, then there exist at least two
fourtuples $y_{i}, (i = 1,...,4)$  -- different from
$x_{i}, (i = 1,...,4)$,  and such that
$x_{1} + x_{2} - x_{3} - x_{4}$ and
$y_{1} + y_{2} - y_{3} - y_{4}$ are linearly equivalent
on $B$. The last is true (for general $x_{i}$'s)
{\it iff}  $\{ x_{1},x_{2} \} = \{ y_{1},y_{2} \}$
and $\{ x_{3},x_{4} \} = \{ y_{3},y_{4} \}$.  This
coincidence of non-ordered 2-point sets
implies the coincidence of the corresponding degenerate
conics -- contradiction. Therefore,
$Supp({\Theta}) = {\psi}({\cal C}_{+})$,
{\bf q.e.d.} Therefore:

\smallskip

{\bf (3).}
{\sl
The theta divisor $\Theta$ of $J(X)$ is
a translate of the Abel-Jacobi image of the
$6$-dimensional family ${\cal C}^{0}_{2,4} (X)$
of non-isolated minimal sections of the system ${\cal S}[2]$.
}

\smallskip

Just as in the proof of the contradiction from above,
it follows:

\smallskip

{\bf (4).}
{\sl
${\Theta}$ is a translate of
${\Phi}^{+,+}({\cal F}_{-} \times {\cal F}_{-})$,
where
${\Phi}^{+,+} : (C_{1},C_{2}) \mapsto {\Phi}([C_{1} + C_{2}])$.
}

\bigskip
{\bf (6.4.3)}
{\it Parameterization of $Sing({\Theta}).$}

\smallskip

The well-known from the theory of jacobians of curves set
$Sing({\Theta}) \cong Sing({\Theta}(B))$
can be expressed also in terms of the
Prym variety $P(\tilde{\Delta}, {\Delta}) \cong J(B)$.

The first which has to be noted is that the planes
of the conics $q(f_{x}), x \in {\bf P}^{2}$ have a common point
${\xi}_{0} \in B$. Indeed, any such  $q(f_{x})$ is
a 6-secant of the curve $B$, and $deg(B) = 7$.
Therefore the plane $<q(f_{x})> \subset {\bf P}^{3}$
intersects $B$ in an additional point, different
from the $6$ common points of $B$ and $q(f_{x})$,
and we have to see that this point does not depend on $x$.

Let $l$ be a general
line in ${\bf P}^{2}$.
Then the map
${\phi}_{l} : l \rightarrow B$,
$x \mapsto B \cap Span(q(f_{x})) -  B \cap q(f_{x})$
is regular, and the curve $B$ is non-rational. Therefore
the image  ${\phi}_{l}(l) \subset B$, of the rational curve $l$,
is a point.  Clearly, this
point ${\xi}_{0}$ does not depend on the choice of $l$
(since any two lines in ${\bf P}^{2}$ intersect each other).

Now,  it is easy to see that
$K_{B}$ = (the hyperplane section of $B$) + ${\xi}_{0}$.

In particular $B$ is a projection of the canonical curve
$B_{can}$ from ${\xi}_{0}$;  and the the Riemann-Roch
theorem implies that the singularities
of ${\Theta}(B)$ correspond to the non-trivial linear
systems defined by some family of fourtuples of
coplanar points on $B$. It is not hard to see, that in
the terms of the attached Prym variety,
these singularities of $\Theta$ can be described
by the family of
degenerations of $C \in {\cal C}^{0}_{2,4}$ to quasi-sections
of type $C = C^{0}_{2,2} + f_{x}$, where
$C^{0}_{2,2}$ is a rational curve on $X$ of bidegree
$(2,2)$. Since all the $(0,2)$ conics on $X$ are
rationally equivalent, we obtain:

\smallskip

{\bf (5).}
{\sl
The Abel-Jacobi image ${\Phi}({\cal C}^{0}_{2,2})(X)$
of the family of rational curves on $X$ of bidegree
$(2,2)$ coincides with the non-trivial
$1$-dimensional component of $Sing({\Theta})$.
}

\smallskip

By the theory of jacobians of curves,
this component is naturally
isomorphic to the curve $\tilde{\Delta}$
(see e.g.[ACGH, Ch.VI, Ex.F]).

\smallskip

We shall see that similar degenerations describe
components of $Sing({\Theta})$  for the bidegree
$(2,2)$ threefold and for the quartic double solid
(see (6.5) and (6.6)).

\bigskip
\bigskip

{\bf (6.5)}
{\sc
The bidegree $(2,2)$ threefold $X$} (see [Ve], [I]).

\smallskip

{\bf (6.5.1)}
{\it The two conic bundle structures on $X$.}

\smallskip

Let $W \subset {\bf P}^{8}$ be the Segre
fourfold
${\bf P}^{2} \times {\bf P}^{2}$, and let
$X$ be an intersection of $W$ with a general quadric,
i.e. $X$ is a {\it bidegree $(2,2)$ threefold}.

Let $p$ and $q$ be the two standard projections
from $W$ to ${\bf P}^{2}$,
(resp. -- from $X$ to ${\bf P}^2$).
Clearly, $p$ and $q$ define conic bundle structures on $X$.

Let
$l = [p^{*}({\cal O}(1)])$ and
$h = [q^{*}({\cal O}(1)])$ be the generators of
${\bf Pic}(W)$ (resp. - of ${\bf Pic}(X)$).

Call the 1-cycle $C$ on $X$ a {\it bidegree $(m,n)$-cycle},
if $C$ has degree $m$ with respect to $l$,  and
degree $n$  --  w.r. to $h$.

\smallskip

{\bf (6.5.2)}
{\it The families ${\cal C}_{+}$ and ${\cal C}_{-}$
for $p.$}

\smallskip

Fix the projection, say $p$. Then
$p : X \rightarrow {\bf P}^{2}$ is a standard conic
bundle, and the discriminant $\Delta$ is a smooth
general plane sextic.

Therefore, the jacobian $J(X)$ is a $9$-dimensional Prym
variety. The two marked families ${\cal C}_{-}$ and
${\cal C}_{+}$ are the families of isolated and
non-isolated minimal sections of the system
${\cal S}[3]$ of conic bundle surfaces,  over the
${\bf P}^{9}$-space of plane cubics.
Theorem(4.4) tells that the Abel-Jacobi
image of one of these two families is a copy of $\Theta$.
It is proven in [I] that the family which parameterizes
the theta divisor is  ${\cal C}_{+}$.  Since the
proof of this coincidence and the description of the
families ${\cal C}_{-}$ and ${\cal C}_{+}$ are purely
technical, and do not differ substantially from the
study in (6.1) - (6.4), we shall only sketch in brief the
results (see [I]):

\smallskip

Let ${\cal C}^{g}_{m,n}$ be the family,
the general element of which is a smooth
connected curve
$C \subset X$ of genus $g$ and of bidegree $(m,n)$;
e.g.,
${\cal C}^{0}_{0,1} \cong \tilde{\Delta}$;
${\cal C}^{0}_{0,3} = \oslash$, etc.
Then the following takes place:

\smallskip

  {\bf (1).}
{\sl
Let
${\cal C}_{+}$ be the $10$-dimensional family of
non-isolated minimal sections, and let
${\cal C}_{-}$ be the $9$-dimensional family of isolated
minimal sections of the system of conic
bundle surfaces ${\cal S}[3]$. Then:
${\cal C}_{+} = {\cal C}^{1}_{3,7}$,
${\cal C}_{-} = {\cal C}^{1}_{3,6}$.
}

\smallskip

  {\bf (2).}Parameterization of ${\Theta}$:

{\sl
{\bf (i).}
$Supp({\Theta}) = C_{+}$. Therefore
${\Phi}_{+}({\cal C}^{1}_{3,7})$ is a copy of $\Theta$;

{\bf (ii).}
${\Phi}_{-}({\cal C}^{1}_{3,6})$ coincides with $J(X)$.
}

In particular  (1), (2), and Theorem (5.8) imply:

\smallskip

(after identifying the points
$z \in Z_{+} \cup Z_{-}$
and the corresponding sheaves
${\cal L} \in {\Theta} \cup P^{-})$:

  {\bf (3).} The fibers of the Abel-Jacobi maps
${\Phi}_{+}$ and ${\Phi}_{-}$:

{\sl
The general fiber $({\Phi}_{+})^{-1}({\cal L})$
of ${\Phi}_{+}$ is an elliptic fibration
over the rational base
$\mid {\cal L} \mid \cong {\bf P}^{1}$.
Let $L \in \mid {\cal L} \mid$ be general. Then
the fiber ${\psi}^{-1}(L)$ is isomorphic
to the plane cubic $C_{0}(L)$ = (the only plane cubic
which passes through the $18$ points
$p(L) \cap {\Delta}$).

The surjective map
${\Phi}_{-} : {\cal C}^{1}_{3,6} \rightarrow J(X)$
is generically finite of degree $2$. If
${\cal L} \in P^{-}$ is general, and $L$ is the
unique effective divisor in
$\mid {\cal L} \mid$,
then
$({\Phi}_{-})^{-1}({\cal L}) \cong
{\psi}^{-1}(L)$  = the two minimal sections
of the ruled surface $S(L).$
}

\smallskip

  {\bf (4).} Parameterization of ${Sing}^{st}{\Theta}$ via
degenerate sections:

{\sl
Let ${\cal C}^{1}_{3,3}$ be the ($6$-dimensional) family
of elliptic curves of bidegree $(3,3)$ on $X.$
Then the Abel-Jacobi image
$Z = {\Phi}({\cal C}^{1}_{3,3})$ is a $3$-dimensional
component of stable singularities of $\Theta$.
}

\bigskip

An elliptic sextic $C_{3,3}$ of bidegree (3,3)
can be treated also as a component of
a (connected) quasi-section
$C_{3,7} = C_{3,3} + \mbox{ two fibers of } p$.

Let
${\cal C}_{sing} \subset {\cal C}^{1}_{3,7}$
be the set of all such quasi-sections.
Since all the fibers of $p$ are rationally equivalent,
the subset $Z \subset J(X)$ can be also treated
as a translate of the set
${\Phi}_{+}({\cal C}_{sing})$.

We shall see in (6.6.5) that similar quasi-sections describe
a singular component of $\Theta$ for the nodal quartic
double solid. The difference is that the theta divisor
of the nodal q.d.s. is described by isolated sections.
This last result is mentioned implicitly in [C2],
in a different context.

\bigskip
\bigskip

  {\bf (6.6)}
{\sc
The nodal quartic double solid.
}

\smallskip

{\bf (6.6.1)}
By definition, a quartic double solid is a double covering
${\rho} : X \rightarrow {\bf P}^{3}$ branched along a quartic
surface $B \subset {\bf P}^{3}$.

The quartic double solid (see [C1], [C2],
[W1], [Vo], [T], [De])
is the most popular example of a Fano threefold, together with
the cubic and the intersection of three quadrics. In particular,
the nodal q. d. s. has a birational structure of a conic bundle
(see e.g. (6.6.3)).

Here we shall use the theory from sections $1-5$, in order to
find additional information about the known parameterization
(see (A)) of $\Theta$ for the general q.d.s.
(see Corollary (6.6.6)(*)).

\bigskip

{\bf (6.6.2).}
{\it Summary of some known results}
(see [T], [C2], [De], [Vo]).

\smallskip

{\sl
Let the the branch locus $B$ has $0 \le {\delta} \le 7$
simple nodes, which impose independent conditions on the
${\bf P}^{9}$-space of quadrics in ${\bf P}^{3}$.
Let $R \cong B$ be the ramification divisor, and let
${\Sigma} = Sing(R) = Sing(X) \cong Sing(B)$ be the
set of nodes of $X$. Let ${\cal R}_{\Sigma}$ be the
$(12 - {\delta})$-dimensional family of Reye sextics
(sextics of genus $3$) on $X$,  which pass through
all the points of $\Sigma$
(and also the proper
preimage of ${\cal R}_{\Sigma}$ on the desingularization
$\tilde{X}$ of $X$).
Let $(J,{\Theta})$ be the $(10 - {\delta})$-dimensional
p.p. intermediate jacobian of $\tilde{X}$. Then:

 {\bf (A).}  The Abel-Jacobi map
${\Phi}_{\cal R} : {\cal R}_{\Sigma} \rightarrow J$
sends ${\cal R}_{\Sigma}$ onto a copy of $\Theta$,
and the connected components of the general fiber
of ${\Phi}_{\cal R}$ are isomorphic to the projective
space ${\bf P}^{3}$.

  {\bf (B).} Let ${\delta} \le 5$, and let ${\cal E}_{\Sigma}$
be the $(8 - {\delta})$-dimensional family of elliptic
quartics on $X$ which pass through all the points of $\Sigma$
(and also the proper preimage of
${\cal E}_{\Sigma}$ on $\tilde{X}$). Then:

The Abel-Jacobi map
${\Phi}_{\cal E} : {\cal E}_{\Sigma} \rightarrow J$
sends ${\cal E}$ onto a $(5 - {\delta})$-dimensional
component $Z$ of $Sing \ {\Theta}$.
Let $\tilde{B} \subset {\bf P}^{9-{\delta}} \cong
{\bf P}( \mbox{ Tang.space of } J \mbox{ in } o)$
be the image of $S$ by the non-complete system
$\mid {\cal O} (2 - Sing(B)) \mid$ of quadrics through $Sing(B)$.
Then:

The general
$z \in Z$ is a quadratic singularity of $\Theta$;
the tangent cone $Cone_{z}({\Theta})$ has rank $5$, and
$Cone_{z}$ passes through $\tilde{B}$. Moreover, any quadric
of rank $5$ through $\tilde{B}$ arises from some $z \in Z$,
and the intersection of all such cones
$Cone_{z}$ coincides with $\tilde{B}$ (= the Torelli
theorem for $\tilde{X}$).
}

\bigskip

{\it Comments.}

Let $X$ be smooth, i.e. ${\delta} = 0$.
It follows from [T] and [C2] that the connected
components of the general fiber of ${\Phi}_{\cal R}$
are ${\bf P}^{3}$-spaces, and any
such a component coincides with a complete linear
system of Reye sextics on $X$ which lie on a fixed
$K3$-surface
$S \in \mid {\cal O}_{X} (2) \mid \cong
{\bf P} {\bf ( } {\rho}^{*} H^{0}({ {\cal O}_{P^{3}} } (2) )
+ {\bf C}.R {\bf )}$ (i.e. $S$ is a quadratic section of $X$).

It is also not hard to see that the family of quadratic sections
of $X$ containing Reye sextics, form a codimension 1 subvariety
of the complete linear system ${\bf P}^{9}$ of all quadratic
sections of $X$.

\bigskip

  {\bf (6.6.3).}
{\it The conic bundle structure on the nodal q.d.s.}

\smallskip

Let $S$ has a simple node $o$.  Denote by $o$ also
the node of $X$ -- ``above'' $o$.
Let
$\tilde{B} \subset \tilde{\bf P} \subset {\bf P}^{8}$
be the image of $B \subset {\bf P}^{3}$ by
the system of quadrics through $o$, and let
$\tilde{\rho} : \tilde{X} \rightarrow \tilde{\bf P}$
be the induced double covering branched along
$\tilde{B}$. (The threefold $\tilde{\bf P}$ is a projection
of the Veronese image ${\bf P}^{3}_{8} \subset {\bf P}^{9}$
of ${\bf P}^{3}$, through the image of $o$. In particular,
$\tilde{\bf P}$ contains a plane ${\bf P}^{2}_{0}$, and
the inverse map
${\sigma} : \tilde{\bf P} \rightarrow {\bf P}^{3}$
is a blow-down of ${\bf P}^{2}_{0}$ to $o$. The
restriction ${\sigma} : \tilde{B} \rightarrow B$
is a blow-down of a smooth conic $q_{o} \subset \tilde{B}$
to the node $o$.)

The threefold
$\tilde{\bf P} \cong
{\bf P}_{P^{2}} ( {\cal O} \oplus {\cal O} (1) )$
has a natural projection $p_{o}$ to

${\bf P}^{2} = \{ \mbox{ the lines } l \mbox{ in } {\bf P}^{3}
\mbox{ through } o \}$,
and ${\bf P}^{2}_{0}$ is the exceptional section
of the projectivized bundle $\tilde{\bf P}$.
The general fiber $p^{-1}(l)$ of the composition
$p = p_{o} \circ \tilde{\rho} :
\tilde{X} \rightarrow {\bf P}^{2}$
is a smooth conic
$q(l) = p^{-1}(l) \cong$ (the
desingularization of ${\rho}^{-1}(l)$ in $o$).

The restriction
${p_{o}\mid}_{\tilde{B}}:\tilde{B} \rightarrow {\bf P}^{2}$
desingularizes the projection from the quartic $B$ through
the node $o = Sing(B)$.
Therefore, ${p_{o}\mid}_{\tilde{B}}$ is a double
covering branched along a smooth plane sextic $\Delta$, and
the conic $q_{o}$ is totally tangent to $\Delta$.
Clearly, the fiber $p^{-1}(x)$ is singular for any
$x \in {\Delta}$, and the natural Abel-Jacobi map
$\tilde{\Delta} \rightarrow J = J(\tilde{X})$ induces an
isomorphism of p.p.a.v.
$P(\tilde{\Delta},{\Delta}) \cong J$
(see [B1]).

\smallskip

{\bf (6.6.4)}
{\it The families ${\cal C}_{+}$ and ${\cal C}_{-}$.}

\smallskip

It is not hard to find the families ${\cal C}_{+}$ and
${\cal C}_{-}$ for $p$.  Since this description
does not differ substantially from the general one, we
shall state it in a brief:

Since $\tilde{\bf P} \subset {\bf P}^{8}$, the degree
map
$deg : \{  \mbox{ subschemes of } \tilde{\bf P} \}
\rightarrow {\bf Z}$
is well defined. In particular,
$deg(\tilde{\bf P}) = deg({\bf P}^{3}_{8}) - 1 = 7$.

Let $Z \subset \tilde{X}$ be a subscheme of $\tilde{X}$.
Define $deg(Z) := deg({\tilde{\rho}}_{*}(Z)).$

{\it Example.}
Let  $l \subset {\bf P}^{3}$ be a line through $o$, let
$x = [l] \in {\bf P}^{2}$ be the point representing $l$,
and let $q(x)$ be the ``conic'' $q(x) = p^{-1}(x)$.
Then $deg(p^{-1}(x)) = 2$.  Indeed,
${\tilde{\rho}}_{*}(q(x)) = 2l'$ where
$l' = p_{o}^{-1}(x) \subset \tilde{\bf P}$
is the line in ${\bf P}^{8}$ which represents
the ``bundle-fiber'' in $\tilde{\bf P}$ over $[l]$.  Note also
that $l'$ is the proper preimage of the line
$l \subset {\bf P}^{3}$ under the blow-down
${\sigma}_{o} : \tilde{\bf P} \rightarrow {\bf P}^{3}$.

\smallskip
{\bf Proposition.}
{\sl
Let ${\cal S}[3]$ be the canonical family of conic
bundle surfaces over plane cubics, induced by $p$,
and let ${\cal C}_{+}$ and ${\cal C}_{-}$ be the
$10$-dimensional family of non-isolated minimal
sections and the $9$-dimensional family of
isolated minimal sections of the system ${\cal S}[3]$.
Then:
}
\smallskip

  {\bf (1).}
{\sl
${\cal C}_{+}$ is a component ${\cal C}^{1}_{10}$
of the family of elliptic curves of degree $10$
on $\tilde{X}$.
}

\smallskip

  {\bf (2).}
{\sl
${\cal C}_{-}$ is a component ${\cal C}^{1}_{9}$
of the family of elliptic curves of degree $9$
on $\tilde{X}$.
}

\smallskip

  {\bf (3).}
{\sl
${\psi}({\cal C}_{-}) = Supp({\Theta})$.
Therefore
${\Phi}_{-}({\cal C}^{1}_{9})$ is a copy of $\Theta$;
${\Phi}_{+}({\cal C}^{1}_{10}) = J$.
}

\smallskip

  {\bf (4).}
{\sl
The general fiber ${\Phi}_{-}^{-1} ({\cal L})$ is
isomorphic to a disjoint union of two lines
${\bf P}'$ and ${\bf P}''$; and the elements of
${\bf P}'$ and ${\bf P}''$ are fibers of elliptic
fibrations on two fixed $K3$-surfaces
$S'$ and $S''$ on $\tilde{X}.$
}

\smallskip

{\bf Proof.}
The proof of (1) and (2) is standard.

\smallskip

{\it Proof of (3).}
The general element $C \in {\cal C}_{-}$ lies
in a unique $S$ of the $9$-dimensional
system $\mid {\cal O}_{\tilde X} (1) \mid$ =
$\mid \tilde{\rho}^{*} {\cal O}_{\tilde{\bf P}} (1) \mid$.
Since $S = S(C)$ is a $K3$-surface, and $C$ is an elliptic
curve on $S$, the base of the linear system
$\mid C \mid$ on $S$ is isomorphic to ${\bf P}^{1}$.

The map
${\psi} : \mid C \mid \rightarrow
Supp({\Theta}) \cup Supp(P^{-})$ sends
the elliptic fibration $\mid C \mid$
onto a non-trivial linear system ${\cal L}$.
Therefore, ${\psi}(C) \in Supp({\Theta})$
{\sl q.e.d.}

\smallskip

{\it Proof of (4).}
It follows from (3), and from Theorem (5.8),
that ${\psi}$ is generically of degree $2$.  Moreover,
if ${\cal L} \in \Theta$  and $L \in \mid {\cal L} \mid$
are general,  then ${\psi}^{-1}(L) = \{ C',C'' \}$, where
$C'$ and $C''$ are the two minimal sections (zero-sections)
of the defined by $L$ ruled surface $S(L)$.
It is not hard to see that the general $C'$ (respectively
$C''$) lies on a unique $K3$-surface $S'$ (resp. $S''$),
and defines there an elliptic fibration
$\{ C'(t) \}$ (resp. $\{ C''(t) \}$).
Now, (4) follows from the preceding and from Theorem (5.8).
{\bf q.e.d.}

\smallskip

  {\bf (6.6.5).} The component $Z \subset Sing({\Theta})$:

{\bf Proposition.}
{\sl
Let $X$ has a simple node $o$.
Let ${\cal E}$ be the $8$-dimensional family of elliptic
quartics on $X$, and let ${\cal E}_{0}$ be the
$7$-dimensional subfamily of elliptic quartics on $X$
through $o$. Then the proper preimage of ${\cal E}_{0}$
on $\tilde{X}$ is a component ${\cal C}^{1}_{7}$
of elliptic curves of degree $7$ on $\tilde{X}$.
The conic bundle projection
$p : \tilde{X} \rightarrow {\bf P}^{2}$
maps the general $C \in {\cal C}^{1}_{7}$
isomorphically onto a plane cubic $p(C)$.
Therefore the Abel-Jacobi map
${\Phi} : {\cal C}^{1}_{7} \rightarrow J = J(\tilde{X})$
sends ${\cal C}^{1}_{7}$ onto a $4$-dimensional
component $Z$ of $Sing({\Theta})$.
}

\smallskip

{\bf Proof. }
Since the proof of (8) can be obtained by a standard
degeneration argument from (B), we shall only sketch
it in brief.

The elements of ${\cal C}^{1}_{7}$ form a family of
isolated minimal sections on surfaces of the system
${\cal S}[3]$.  However this family
is not globally defined on the system
${\cal S}[3]$ (i.e. the general surface of this system
does not admit such a section).
The general element
$C \in {\cal C}^{1}_{7}$ can be completed
(by  $\#(p(C) \cap {\Delta})$ ways) to a connected
quasi-section $(C + \mbox{ fiber}) \in {\cal C}^{1}_{9}$.
Since all the fibers of $p$ are rationally equivalent,
it follows that $Z$ can be treated as a translate
of the set ${\Phi}_{-}({\cal C}_{sing})$ where
${\cal C}_{sing} =
\{ C + \mbox{ fiber } \} \cap {\cal C}^{1}_{9}$.

It is not hard to see that the points $z$ of the
$4$-dimensional set $Z := {\Phi}_{{\cal C}^{1}_{7}}$
describe stable singularities
${\cal L} = {\cal L}(z)$ of ${\Theta}$
regarded as a theta divisor on a Prym variety.

Indeed, the general $C \in {\cal C}^{1}_{7}$ lies
on a ${\bf P}^{2}$-system of $K3$-surfaces
$S(u,v)$ of the system
$\mid {\cal O}_{\tilde{X}} (1) \mid$,
and any embedding $C \subset S(u,v)$
includes the curve $C$ in an elliptic pencil
$C(u,v;t)$ on $S(u,v)$.
Now, the map
${\psi} : \{ C(u,v;t) \} \rightarrow
Supp({\Theta}) \cup Supp(P^{-})$
defines a $3$-dimensional linear system
$\{ L(u,v;t) \}$ on $\tilde{\Delta}$.
{bf q.e.d.}

\bigskip

{\bf (6.6.6).}
{\bf (c).}
{\it From nodal q.d.s to general q.d.s.}

\smallskip

The parameterization of $\Theta$ for the general q.d.s.
can be obtained from the parameterization of $\Theta$ for
the nodal q.d.s.
The consideration offered here is based
on the Clemens's suggestion to prove the parameterization
of $\Theta$ of the general q.d.s by using degeneration to
a nodal q.d.s. ([see [C2, p.98]).

Let $X_{t}$ be a general
Lefschetz pencil of quartic double solids, s.t. $X_{0}$
is a nodal q.d.s.  One can consider $\{ X_{t} \}$ as a
pencil of hyperplane sections
in a fourfold ${\pi}:W \rightarrow {\bf P}^4$,
which is a double covering of
${\bf P}^4$ branched along a smooth quartic 3-fold.

Let
$\tilde{X_{0}} \rightarrow {\bf P}^2$ be the induced conic
bundle, and let ${\cal C}_{-}$ be the
familiy of isolated  minimal sections on
$\tilde{X_{0}}$.

Define the degree of the effective 1-cycle $C \subset W$
as the degree of the effective 1-cycle ${\pi}_{*}(C)$.

It is not hard to see that image on $X_{0}$
of the family ${\cal C}_{-}$ coincides with the
9-dimensional family ${{\cal C}^{3}_{6}}[3]$ of
curves $C \subset X_{0}$ s.t. $deg(C) = 6$,
$p_{a}(C) = 3$, and $C$ has a triple point in the node
$o$ of $X_{0}$. The Abel-Jacobi image of this family coincides
with the Abel-Jacobi image of the 12-dimensional family
${\cal R}(X_{0})$ of Reye sextics on $X_{0}$, since any Reye
sextic on $X_{0}$ can be included in a rational family
($\cong {\bf P}^3$) of Reye sextics containing an element
of ${{\cal C}^{3}_{6}}[3]$ (see e.g. [C2]).  Indeed,
rationally equivalent 1-cycles have the same Abel-Jacobi image.
The rest repeats the arguments from [C2].

However, the obtained description (4) of the fiber
of ${\Phi}_{-}$ adds a new information to the known results
(A) and (B). More concretely, let $X = X_{t}$ be a general
q.d.s., and let ${\cal R}$ be the family of Reye sextics
on $X$.  According to (A), any
connected component of the general fiber of
${\Phi}: {\cal R} \rightarrow {\Theta}$ is isomorphic to
${\bf P}^{3}$. Now, (4) implies that these components are
exactly {\it two}.

In just the same way one can see that the image on $X_{0}$ of
family ${\cal C}_{+}$ is a subset of the 14-dimensional
family ${\cal D}(X_{0})$ of septics of genus 4 on $X_{0}$.
This, in addition, proves the existence of such curves
on the general $X_{t}$. Now, the same argument as above imply
that the Abel-Jacobi map sends the corresponding family
${\cal D}(X_{t})$ surjectively onto the intermediate jacobian
$J(X_{t})$.

\smallskip

{\bf (6.6.7).}
{\bf Corollary.}
{\sl
Let $X$ be a general smoooth quartic double solid,
and let $(J(X), {\Theta})$ be the principally polarized
intermediate jacobian of $X$. Then

{\bf (i).}
The Abel-Jacobi map
${\Phi}_{\cal R}$
sends the 12-dimensional family
${\cal R}$
of Reye sextics (sextics of genus 3) on $X$ onto a copy of
${\Theta}$.  Moreover, the general fiber of ${\Phi}_{\cal R}$
has two connected components -- each isomorphic to the
projective space ${\bf P}^{3}$.

{\bf (ii).} The Abel-Jacobi map ${\Phi}_{\cal D}$ sends
the 14-dimensional family ${\cal D}$ of septics of genus 4
on $X$ surjectively onto the intermediate jacobian $J(X)$.
}

\bigskip
\bigskip

\centerline{\sc References}

\bigskip

[ACGH] E.Arbarello, M.Cornalba, P.A.Griffiths, J.Harris,
{\it Geometry of Algebraic Curves, Vol.I.}
Springer-Verlag, New York Inc. (1985).

\smallskip

[B1]    A.Beauville, {\it Vari\'et\'es de Prym et jacobienne
interm\'ediaires.}
Ann.de l'AENS, 4 ser.,10 (1977), 149-196.

\smallskip

[B2]    A.Beauville, {\it Les singularit\'es du diviseur theta
de la jacobienne interm\'ediaire de l'hypersurface cubique
dans ${\bf P}^{4}$.}
Lecture Notes in Math., Vol. 947 (1982), 190-208.

\smallskip

[B3]    A.Beauville, {\it Sous-vari\'et'es speciales des
vari\'et\'es de Prym.}
Compositio Math., 45 (1982), 357-383.

\smallskip

[BM]    S.Bloch, J.P.Murre, {\it On the Chow group of
certain types of Fano threefolds.}
Compositio Math., 39 (1979), 47-105.

\smallskip

[C1]    H.Clemens, {\it Double Solids.}
Advances in Math. 47 (1983), 107-230.

\smallskip

[C2]    H.Clemens, {\it The Quartic Double Solid Revisited.}
Proc. Symp. in Pure Math., Vol. 53 (1991), 89-101.

\smallskip

[CG]    H.Clemens, P.Griffiths, {\it The intermediate jacobian
of the cubic threefold.}
Annals of Math., 95 (1972), 281-356.

\smallskip

[De]    O.Debarre, {\it Sur le theoreme de Torelli pour les
solides doubles quartiques.}
Compositio Math., 73 (1990), 161-187.

\smallskip

[Do]    R.Donagi, {\it Group law on the intersection of
two quadrics.}
Ann. Scuola Norm. Super. Pisa: ser.4, vol. 7:2 (1980),
217-239.

\smallskip

[H]     R.Hartshorne, {\it Algebraic Geometry.}
Springer-Verlag (1977).

\smallskip

[Isk]   V.A.Iskovskikh, {\it On the rationality problem for
conic bundles.}
Duke Math.Journal, 54:2 (1987), 271-294.

\smallskip

[I]    A.Iliev, {\it The theta divisor of the bidegree (2,2)
threefold in ${\bf P}^{2} \times {\bf P}^{2}$.}
preprint (1994).

\smallskip

[LN]    H.Lange, M.S.Narasimhan, {\it Maximal Subbundles of
Rank Two Vector Bundles on Curves.}
Math. Ann., 266 (1983), 55-72.

\smallskip

[M]     S.Mori, {\it Threefolds whose canonical bundles are
not numerically effective.}
Ann of Math., 116 (1982), 133-176.

\smallskip

[Miy]   Y.Miyaoka, {\it On the Kodaira dimension of minimal
threefolds.}
Math. Ann., 281:No.2 (1988), 325-332.

\smallskip

[MM]    S.Mori, S.Mukai, {\it On Fano 3-folds with
$B_{2} \ge 2$.}
Adv. St. in Pure Math., Vol. I - Algebraic Varieties and
Analytic Varieties. Kinokuniya Comp. (1983), 101-129.

\smallskip

[S]     V.G.Sarkisov, {\it Birational automorphisms of
conical fibrations.}
Izv. Akad. Nauk SSSR: Ser.Mat., 44 (1980), 918-945
= Math. USSR-Izv. 17 (1981).

\smallskip

[Se]    W.Seiler, {\it Deformations of ruled surfaces}.
J. Reine Angew. Math., 426 (1992), 203-219.

\smallskip

[T]     A.Tikhomirov, {\it The Abel-Jacobi map of sextics of genus
$3$ on double spaces of ${\bf P}^{3}$ of index two}.
Soviet Math. Dokl., Vol. 33:1 (1986), 204-206.

\smallskip

[Tju]   A.Tjurin, {\it The Middle Jacobian of Three-Dimensional
Varieties}. Journal of Soviet Mathematics,
Vol.13, No.6 (1980), 707-744.

\smallskip

[Ve]    A.Verra, {\it The Prym map has degree two on plane sextics}.
preprint (1991).

\smallskip

[Vo]    C.Voisin, {\it Sur la jacobienne interm\'ediaire du double
solide d'indice deux}.
Duke Math. Journal, 57:2 (1988), 629-646.

\smallskip

[W1]    G.E.Welters, {\it Abel-Jacobi isogenies for certain types
of Fano threefolds}.
Math. Centre Tracts 141, Math. Centrum  Amsterdam (1981).

\smallskip

[W2]    G.E.Welters, {\it A theorem of Gieseker-Petri type
for Prym varieties}.
Ann. Sci. E.N.S., 4 ser., 18 (1985), 671-683.

\smallskip

[Z]     A.A.Zagorskii, {\it three-dimensional conic bundles.}
Math. Notes Akad. Sci. USSR, 21(1977), 420-427.

\smallskip

[ZR]    Ziv Ran, {\it Hodge theory and the Hilbert scheme.}
J. Diff. Geometry, 37:1 (1993), 191-198.

\bigskip
\bigskip
\bigskip
\bigskip

Atanas Iliev

Institute of Mathematics, Bulgarian Academy of Sciences,

Acad.G.Bonchev Str., bl.8,   1113  Sofia,  Bulgaria.

E-mail address:  algebra@bgearn.bitnet

\end{document}